\documentclass{article} 
\usepackage{iclr2023_conference,times}

\usepackage[utf8]{inputenc} 
\usepackage[T1]{fontenc}    
\usepackage{hyperref}       
\usepackage{url}            
\usepackage{booktabs}       
\usepackage{amsfonts}       
\usepackage{microtype}      
\usepackage{color,xcolor}         

\usepackage{makecell}
\usepackage{soul} 
\definecolor{myyellow}{rgb}{1,1,0.75}
\sethlcolor{myyellow}

\usepackage{mathtools}   
\usepackage{algorithm}
\usepackage{algorithmic}

\definecolor{lgray}{rgb}{0.95,0.95,0.95}
\definecolor{yel}{rgb}{1,0.98,0.92}
\definecolor{mydarkblue}{rgb}{0,0.2,0.9}
\newcommand\bl[1]{{\color{mydarkblue}#1}}
\newcommand\ff{\bl{f}}
\newcommand\ffi{\bl{f_i}}
\newcommand\ffc{\bl{f^*}}
\definecolor{mydarkred}{rgb}{0.8,0.0,0.0}
\newcommand\dr[1]{{\color{mydarkred}#1}}
\newcommand\g{\dg{g}}
\newcommand\gc{\dg{g^*}}
\definecolor{mydarkgreen}{rgb}{0,0.55,0}
\newcommand\dg[1]{{\color{mydarkgreen}#1}}
\newcommand\h{\dr{h}}
\newcommand\hc{\dr{h^*}}
\newcommand\hi{\dr{h_i}}
\newcommand\hic{\dr{h_i^*}}

\usepackage{xspace}
\usepackage{scalefnt}
\definecolor{myorange}{RGB}{255,100,0}
\newcommand{\algn}[1]{{\sf\color{mydarkred}\scalefont{0.90}{#1}}\xspace}

\usepackage{amsmath, amsfonts,amssymb}
\usepackage{cite}
\DeclareMathOperator*{\argmin}{arg\,min}
\DeclareMathOperator*{\minimize}{minimize}
\newcommand{\sqnorm}[1]{\left\| #1 \right\|^2}
\newcommand{\Exp}[1]{\mathbb{E}\!\left[ #1 \right]}
\newcommand{\Rop}{\mathcal{R}}
\newcommand{\Ropp}{\widetilde{\mathcal{R}}}
\newcommand{\oma}{\omega_{\mathrm{ran}}}

\title{RandProx: Primal--Dual Optimization Algorithms with Randomized Proximal Updates}

\author{%
 Laurent Condat \&  Peter Richt\'{a}rik\\
 Visual Computing Center\\
King Abdullah University of Science and Technology (KAUST) \\
Thuwal, Kingdom of Saudi Arabia\\
Contact: see \texttt{https://lcondat.github.io/} 
}

\iclrfinalcopy 
\begin{document}

\maketitle

\begin{abstract}
Proximal splitting algorithms are well suited to solving large-scale nonsmooth optimization problems, in particular those arising in machine learning. We propose a new primal--dual algorithm, in which the dual update is randomized; equivalently, the proximity operator of one of the function in the problem is replaced by a stochastic oracle. For instance, some randomly chosen dual variables, instead of all, are updated at each iteration. Or, the proximity operator of a function is called with some small probability only. A nonsmooth variance-reduction technique is implemented so that the algorithm finds an exact minimizer of the general problem involving smooth and nonsmooth functions, possibly composed with linear operators. We derive linear convergence results in presence of strong convexity; these results are new even in the deterministic case, when our algorithms reverts to the recently proposed Primal--Dual Davis--Yin algorithm. Some randomized algorithms of the literature are also recovered as particular cases (e.g., Point-SAGA). But our randomization technique is general and encompasses many unbiased mechanisms beyond sampling and probabilistic updates, including compression. Since the convergence speed depends on the slowest among the primal and dual contraction mechanisms, the iteration complexity might remain the same when randomness is used. On the other hand, the computation complexity can be significantly reduced. Overall, randomness helps getting faster algorithms. This has long been known for stochastic-gradient-type algorithms, and our work shows that this fully applies in the more general primal--dual setting as well.
\end{abstract}

\section{Introduction}

Optimization problems arise virtually in all quantitative fields, including machine learning, data science, statistics, and many other areas \citep{pal09, sra11, bac12, cev14, pol15, bub15, glo16, cha16, sta16}. In the big data era, they tend to be very high-dimensional, and first-order methods are particularly appropriate to solve them. When a function is smooth, an optimization algorithm typically makes calls to its \textbf{gradient}, whereas for a nonsmooth function, its \textbf{proximity operator} is called instead. Iterative optimization algorithms making use of proximity operators  are called proximal (splitting) algorithms \citep{par14}. Over the past 10 years or so, primal--dual proximal algorithms have been developed and are well suited for a broad class of large-scale optimization problems involving several functions, possibly composed with linear operators \citep{com10,bot14,par14,kom15,bec17,con19,com21,con22}. 

However, in many situations, these deterministic algorithms are too slow, and this is where \textbf{randomized algorithms} come to the rescue; they are variants of the deterministic algorithms with a cheaper iteration complexity, 
obtained by calling a random subset, instead of all, of the operators or updating  a random subset, instead of all, of the variables, at every iteration.
Stochastic Gradient Descent (\textbf{SGD})-type methods \citep{RobbinsMonro:1951,Nemirovski-Juditsky-Lan-Shapiro-2009,Bottou2012,gow20a,gor202,kha20} are a prominent example, with the huge success we all know. They consist in replacing a call to the gradient of a function, which can be itself a sum or expectation of several functions, by a cheaper \textbf{stochastic gradient} estimate. By contrast, replacing the proximity operator of a possibly nonsmooth function by a \textbf{stochastic proximity operator} estimate is a nearly virgin territory. This is an important challenge, because many functions of practical interest have a proximity operator, which is expensive to compute. We can mention the nuclear norm of matrices, which requires singular value decompositions, indicator functions of sets on which it is difficult to project, or optimal transport costs \citep{pey19}.

In this paper, we propose \algn{RandProx} (Algorithm~2), a randomized version of the Primal--Dual Davis--Yin (PDDY) method (Algorithm~1), which is a proximal algorithm proposed recently \citep{sal20} and further analyzed in \citet{con22}. In \algn{RandProx}, one proximity operator that appears in the PDDY algorithm is replaced by a stochastic estimate. \algn{RandProx} is \textbf{variance-reduced} \citep{han19,gor202,gow20a}; that is, through the use of control variates, the random noise is mitigated and eventually vanishes, so that the algorithm converges to an exact solution, just like its deterministic counterpart.  
Algorithms with stochastic errors in the computation of proximity operators have been studied, for instance in \citet{com16}, but the errors are typically assumed to decay or some stepsizes are made decaying along the iterations, with a certain rate.  By contrast, in variance-reduced algorithms such as  \algn{RandProx}, which has fixed stepsizes, error compensation is automatic.

We analyze \algn{RandProx} and prove its linear convergence in the strongly convex setting, with additional results in the convex setting; we leave the nonconvex case, which requires different proof techniques, for future work. We mention relationships between our results and related works in the literature throughout the paper. In special cases,  \algn{RandProx} reduces to  Point-SAGA  \citep{def16}, the Stochastic Decoupling Method \citep{mis19}, ProxSkip, SplitSkip and Scaffnew 
\citep{mis22}, and randomized versions of the PAPC \citep{dro15}, PDHG \citep{cha11a} and ADMM \citep{boy11} algorithms. They are all generalized and unified within our new framework. Thus, \algn{RandProx} paves the way to the design of proximal counterparts of variance-reduced SGD-type algorithms, 
just like Point-SAGA  \citep{def16} is the proximal counterpart of SAGA \citep{def14}.

\section{Problem formulation}

Let $\mathcal{X}$ and $\mathcal{U}$ be finite-dimensional real Hilbert spaces.
We consider the generic convex optimization problem:
\begin{equation}
\mathrm{Find} \ x^\star \in \argmin_{x\in\mathcal{X}}  \Big( \ff(x)+\g(x) + \h(Kx)\Big),\label{eqpb0}
 \end{equation}
where $K:\mathcal{X}\rightarrow \mathcal{U}$ is a nonzero linear operator;
$\ff$ is a convex  $L_{\ff}$-smooth function, for some $L_{\ff}>0$; that is,  its gradient $\nabla \ff$ is $L_{\ff}$-Lipschitz continuous \citep[Definition 1.47]{bau17}; and  $\g:\mathcal{X}\rightarrow \mathbb{R}\cup\{+\infty\}$ and
$\h:\mathcal{U}\rightarrow \mathbb{R}\cup\{+\infty\}$ 
are proper closed convex functions whose proximity operator is easy to compute. 

We will assume strong convexity of some functions: a 
 convex function $\phi$ is said to be $\mu_\phi$-strongly convex, for some  $\mu_\phi\geq 0$, if $\phi-\frac{\mu_\phi}{2}\|\cdot\|^2$ is convex. This covers the case $\mu_\phi=0$, in which $\phi$ is merely convex.

\subsection{Proximity operators and proximal algorithms} 
We recall that for any function $\phi$ and parameter $\gamma>0$, the proximity operator of $\gamma \phi$ is \citep{bau17}:
 $
 \mathrm{prox}_{\gamma \phi}: x\in\mathcal{X} \mapsto \argmin_{x'\in\mathcal{X}}\big(\gamma\phi(x') + \frac{1}{2} \|x'-x\|^2\big)
  $. 
  This operator has a closed form for many functions of practical interest \citep{par14,pus17,ghe18}, see also the website \url{http://proximity-operator.net}. 
   In addition, the Moreau identity holds: 
 \begin{equation*}
\mathrm{prox}_{\gamma \phi^*}(x)=x-\gamma\,\mathrm{prox}_{\phi/\gamma}(x/\gamma),
\end{equation*}
where $\phi^*:x\in\mathcal{X} \mapsto \sup_{x'\in\mathcal{X}} \big(\langle x,x'\rangle -\phi(x')\big)$ denotes the conjugate function of $\phi$ \citep{bau17}. Thus, one can compute the proximity operator of $\phi$ from the one of $\phi^*$, and conversely. \\

Proximal splitting algorithms, such as the forward--backward and the Douglas--Rachford algorithms \citep{bau17}, are well suited to minimizing the sum, $\ff + \g$ or $\g+\h$ in our notation, of two functions.
However, many problems take the form \eqref{eqpb0} with $K\neq \mathrm{Id}$, where $\mathrm{Id}$ denotes the identity, and the proximity operator of $\h\circ K$ is intractable in most cases.  
A classical example is the total variation, widely used in image processing \citep{rud92,cas11,con14,con17} or for regularization on graphs \citep{cou13}, where $\h$ is some variant of the $\ell_1$ norm and $K$ takes 
differences between adjacent values. 
Another example is  
when $\h$ is the indicator function of some nonempty closed convex set $\Omega \subset \mathcal{U}$; that is, $\h(u)=(0$ if $u\in\Omega$, $+\infty$ otherwise$)$, in which case the problem \eqref{eqpb0} can be rewritten as
\begin{equation*}
\mathrm{Find} \ x^\star \in \argmin_{x\in\mathcal{X}} \Big( \ff(x)+\g(x)\Big)\quad\mbox{s.t.}\quad Kx\in \Omega.
 \end{equation*}
 If $\g=0$ and $\Omega=\{b\}$ for some $b\in\mathrm{ran}(K)$, where $\mathrm{ran}$ denotes the range, the problem can be further rewritten as the linearly constrained smooth minimization problem 
\begin{equation*}
\mathrm{Find} \ x^\star \in \argmin_{x\in\mathcal{X}} \, \ff(x) \quad \mbox{s.t.} \quad Kx=b.
 \end{equation*}
This last problem has applications in decentralized optimization, for instance \citep{xin20,kov20,sal22}. Thus, the template problem \eqref{eqpb0} covers a wide range of optimization problems met in machine learning \citep{bac12,pol15}, signal and image processing \citep{com10,cha16}, control \citep{sta16}, and many other fields.  Examples include compressed sensing \citep{can06}, object discovery in computer vision \citep{vo19}, $\ell_1$ trend filtering \citep{kim09}, group lasso \citep{yua06}, square-root lasso \citep{bel11}, Dantzig selector \citep{can07}, and support-vector machines \citep{cor95}.

\subsection{The dual problem, saddle-point reformulation, and optimality conditions}

In order to analyze algorithms solving such problems, 
we introduce the dual problem to \eqref{eqpb0}:
\begin{equation}
\mathrm{Find} \ u^\star \in \argmin_{u\in\mathcal{U}} \, \Big( (\ff+\g)^*(-K^*u)+\hc(u)\Big),\label{eqpbd}
 \end{equation}
 where $K^*:\mathcal{U}\rightarrow\mathcal{X}$ is the adjoint operator of $K$. We can also express the primal and dual problems as a combined saddle-point problem:
\begin{equation}
\mathrm{Find} \ (x^\star,u^\star) \in \arg\min_{x\in\mathcal{X}}\max_{u\in\mathcal{U}} \, \Big( \ff(x)+\g(x)+\langle Kx,u\rangle-\hc(u)\Big).\label{eqpbs}
 \end{equation}
 For these problems to be well-posed, we suppose that 
 there exists $x^\star\in\mathcal{X}$ such that 
\begin{equation}
0\in \nabla \ff(x^\star)+\partial \g(x^\star)+K^*\partial \h(Kx^\star),\label{eqinc0}
\end{equation}
where $\partial (\cdot)$ denotes the subdifferential \citep{bau17}. By Fermat's rule, every $x^\star$ satisfying \eqref{eqinc0} is a solution to \eqref{eqpb0}.
  Equivalently to  \eqref{eqinc0}, we suppose that there exists $(x^\star,u^\star)\in\mathcal{X}\times\mathcal{U}$ such that
  \begin{equation}
 \left\{ \begin{array}{l}
  0\in \nabla \ff(x^\star)+\partial \g(x^\star) +K^* u^{\star}\\ 
0\in -K x^\star+\partial \hc (u^{\star})
\end{array}\right..\label{eqai}
\end{equation}
Every $(x^\star,u^\star)$ satisfying \eqref{eqai} is a primal--dual solution pair; that is, $x^\star$  is a solution to \eqref{eqpb0}, $u^\star$  is a solution to \eqref{eqpbd}, and $(x^\star,u^\star)$ is a solution to \eqref{eqpbs}.

\section{Proposed algorithm: \algn{RandProx}}

There exist several deterministic algorithms for solving the problem \eqref{eqpb0}; see 
\citet{con19} for a recent overview. 
In this work, we focus on the PDDY algorithm (Algorithm~1) \citep{sal20,con22}. In particular, our new algorithm \algn{RandProx} (Algorithm~2) generalizes the PDDY algorithm with a \hl{stochastic estimate of the proximity operator of $\hc$}.

\begin{figure*}[t]
\begin{minipage}{.42\textwidth}
\begin{algorithm}[H]\label{alg:PDDY}
		\caption{
		PDDY algorithm\\ \citep{sal20}}
		\begin{algorithmic}
			\STATE  \textbf{input:} initial points $x^0\in {\cal X}$, $u^0\in \cal U$; \\stepsizes $\gamma>0$, $\tau>0$
			\STATE $v^0\coloneqq K^* u^0$
			\FOR{$t=0, 1, \ldots$}
			\STATE $\hat{x}^{t} \coloneqq  \mathrm{prox}_{\gamma \g}\big(x^t -\gamma \nabla \ff(x^t) - \gamma v^t\big)$%
			\STATE \hl{$u^{t+1}\coloneqq \mathrm{prox}_{\tau \hc} \big(u^t+\tau K \hat{x}^{t}\big)$}\phantom{\large W}
			\STATE $v^{t+1}\coloneqq K^* u^{t+1}$
			\STATE  $x^{t+1} \coloneqq \hat{x}^{t}-\gamma (v^{t+1}-v^t)$
			\ENDFOR
		\end{algorithmic}
	\end{algorithm}\end{minipage}
	\quad
	\begin{minipage}{.54\textwidth}
	\begin{algorithm}[H] \label{alg:RandProx}
		\caption{
		\algn{RandProx}\\ {[new]}}
		\begin{algorithmic}
			\STATE  \textbf{input:} initial points $x^0\in {\cal X}$, $u^0\in {\cal U}$; \\
			stepsizes $\gamma>0$, $\tau>0$; \hl{$\omega\geq 0$} 
			\STATE $v^0\coloneqq K^* u^0$
			\FOR{$t=0, 1, \ldots$}
			\STATE $\hat{x}^{t} \coloneqq  \mathrm{prox}_{\gamma \g}\big(x^t -\gamma \nabla \ff(x^t) - \gamma v^t\big)$
			\STATE \hl{$u^{t+1}\coloneqq u^t + \tfrac{1}{1+\omega} \Rop^t\big(\mathrm{prox}_{\tau \hc} (u^t\!+\!\tau K \hat{x}^{t})-u^t\big)$}\phantom{\Large I}%
			\STATE $v^{t+1}\coloneqq K^* u^{t+1}$ 
			\STATE  $x^{t+1} \coloneqq \hat{x}^{t}-\gamma\,$\hl{$(1+\omega)$}$\,(v^{t+1}-v^t)$ 
			\ENDFOR
		\end{algorithmic}
	\end{algorithm}\end{minipage}\end{figure*}

\subsection{The PDDY algorithm}\label{sec22}

We recall the general convergence result for the PDDY algorithm~\citep[Theorem 2]{con22}: \begin{quote}\emph{If $\gamma\in(0,2/L_{\ff})$, $\tau>0$, $\tau\gamma \|K\|^2\leq 1$, then $(x^{t})_{t\in\mathbb{N}}$ converges to a primal solution $x^\star$ of \eqref{eqpb0} and $(u^{t})_{t\in\mathbb{N}}$ converges to a dual solution $u^\star$ of \eqref{eqpbd}.}\end{quote} 
The PDDY algorithm is similar and closely related to the PD3O algorithm \citep{yan18}, as discussed in \citet{sal20,con22}. It is also an instance (Algorithm 5) of the Asymmetric Forward--Backward Adjoint (AFBA) framework of \citet{lat17}. 
We note that the popular Condat--V\~u algorithm \citep{con13,vu13} can solve the same problem but has more restrictive conditions on $\gamma$ and $\tau$. 

In the PDDY algorithm, the full gradient $\nabla \ff$ can be replaced by a stochastic estimator which is typically cheaper to compute \citep{sal20}. Convergence rates and accelerations of the PDDY algorithm, as well as distributed versions of the algorithm, have been derived   in \citet{con22}. In particular, if $\mu_{\ff}>0$ or $\mu_{\g}>0$, the primal problem  \eqref{eqpb0} is strongly convex. In this case, a varying stepsize strategy  accelerates the algorithm, with a $\mathcal{O}(1/t^2)$ decay of $\|x^{t}-x^\star\|^2$, where $x^\star$ is the unique solution to \eqref{eqpb0}. But strong convexity of the primal problem is not sufficient for the PDDY algorithm to converge linearly, and additional assumptions on $\h$ and $K$ are needed.
We will prove  linear convergence when both the primal and dual problems are strongly convex; this is a natural condition for primal--dual algorithms.

We note that $\h$ is $L_{\h}$-smooth, for some $L_{\h}>0$, if and only if  $\hc$ is $\mu_{\hc}$-strongly convex, for some $\mu_{\hc} >0$, with $\mu_{\hc} = 1/L_{\h}$.  In that case,  the dual problem \eqref{eqpbd} is strongly convex.

\subsection{Randomization mechanism for the proximity operator of $\hc$}

We propose \algn{RandProx} (Algorithm 2), a generalization of the PDDY algorithm (Algorithm 1) with a randomized update of the dual variable $u$. 
Let us formalize the random operations using random variables and stochastic processes. We introduce the underlying probability space $( \mathcal{S}, \mathcal{F}, P)$. Given a real Hilbert space $\mathcal{H}$, an $\mathcal{H}$-valued random variable is a measurable map from $( \mathcal{S}, \mathcal{F})$ to $(\mathcal{H},\mathcal{B})$, where $\mathcal{B}$ is the Borel $\sigma$-algebra of $\mathcal{H}$. Formally, randomizing some steps in the PDDY algorithm amounts to defining $\big((x^t,u^t)\big)_{t\in\mathbb{N}}$ as a stochastic process, with $x^t$ being a $\mathcal{X}$-valued random variable and $u^t$ a $\mathcal{U}$-valued random variable, for every $t\geq 0$. 
We use light notations and write our randomized algorithm \algn{RandProx} 
using 
 stochastic operators $\Rop^t$ on $\mathcal{U}$; 
that is, for every $t\geq 0$ and any $r^t\in\mathcal{U}$, $\Rop^t(r^t)$ is a  $\mathcal{U}$-valued random variable, which can be interpreted as $r^t$ plus `random noise' (formally, $r^t$ is itself a $\mathcal{U}$-valued random variable, but algorithmically, $\Rop^t$ is applied to a particular outcome in $\mathcal{U}$, hence the notation as an operator on  $\mathcal{U}$). 
 To fix the ideas, let us give two examples.\smallskip
 
{\bf Example 1.} The first example is compression \citep{ali17,ali18,Cnat,mis19d,alb20,bez20,con22e}: 
 $\mathcal{U}=\mathbb{R}^d$ for some $d\geq 1$ and $\Rop^t$ is the well known \texttt{rand}-$k$ compressor or sparsifier, with $1\leq k < d$: $\Rop^t$  multiplies $k$ coordinates, chosen uniformly at random, of the vector $r^t$ by $d/k$ and sets the other ones to zero. An application to compressed communication is discussed in Section~\ref{secfl}.\smallskip 
 
{\bf Example 2.} The second example, discussed in Section~\ref{secskip}, is the Bernoulli, or coin flip, operator
\begin{equation}
\Rop^t : r^t \mapsto \begin{cases}\frac{1}{p} r^t & \text{ with probability } p, \\ 0 & \text{ with probability }1-p,\end{cases}\label{eqber}
\end{equation}
for some $p>0$. In that case, with probability $1-p$, the outcome of $\Rop^t(r^t)$ is 0 and $r^t$ does not need to be calculated; in particular, in the \algn{RandProx} algorithm, $\mathrm{prox}_{\tau \hc}$ is not called, and this is why one can expect the iteration complexity of \algn{RandProx} to decrease. Thus, in this example, $\Rop^t(r^t)$ does not really consist of applying the operator $\Rop^t$ to $r^t$; 
\hl{in general,  the notation $\Rop^t(r^t)$ simply denotes a stochastic estimate of $r^t$}.\smallskip

{\bf Example 3.} The third example, discussed in Section~\ref{secsamp}, is sampling, which makes it possible to solve problems involving a sum $\sum_{i=1}^n \hi$ of functions, by calling the proximity operator of only one randomly chosen function $\hi$, instead of all functions, at every iteration. The 
Point-SAGA algorithm \citep{def16} is recovered as a particular case of \algn{RandProx} in this setting.\\

Hereafter, we denote by $\mathcal{F}_t$ the $\sigma$-algebra generated by the collection of $(\mathcal{X}\times\mathcal{U})$-valued random variables $(x^0,u^0),\ldots, (x^t,u^t)$, for every $t\geq 0$. In this work, we consider \textbf{unbiased} random estimates: for every $t\geq 0$, 
\begin{equation*}
\Exp{\Rop^t(r^t)\;|\;\mathcal{F}_t}=r^t,
\end{equation*}
where $\Exp{\cdot}$ denotes the expectation, here conditionally on $\mathcal{F}_t$, and $r^t$ is the random variable $$r^t \coloneqq \mathrm{prox}_{\tau \hc} (u^t+\tau K \hat{x}^{t})-u^t,$$ as defined by \algn{RandProx}. 
Note that our framework is general in that for $t\neq t'$, $\Rop^t$ and $\Rop^{t'}$ need not be independent nor have the same law. In simple words, at every iteration, the randomness  is new but can have a different form and depend on the past, so that the operators $\Rop^t$ can be defined dynamically on the fly in \algn{RandProx}.

We characterize the operators $\Rop^t$ by their \emph{relative variance} 
$\omega\geq 0$ such that, for every $t\geq 0$, 
\begin{equation}
\Exp{\sqnorm{\Rop^t(r^t)-r^t}\;|\;\mathcal{F}_t} \leq \omega \sqnorm{r^t}.\label{eqrelva}
\end{equation}
 This assumption is satisfied by a large class of randomization strategies, which are widely used to define unbiased stochastic gradient estimates. We refer to \citet{bez20}, Table 1 in \citet{saf21}, \citet{zha21}, \citet{sze22} for examples. In the Example 1 above of \texttt{rand}-$k$, $\omega=\frac{d}{k}-1$. In Example 2, $\omega=\frac{1}{p}-1$. In Example 3,   $\omega=n-1$.  The value of $\omega$ is supposed known and is used in the \algn{RandProx} algorithm. Note that $\omega=0$ if and only if $\Rop^t=\mathrm{Id}$, in which case there is no randomness and \algn{RandProx} reverts to the original deterministic PDDY algorithm.

Thus, $\Rop^t(r^t)=r^t+e^t$, with the variance of the error $e^t$ proportional to $\sqnorm{r^t}$. In particular, if $r^t=0$, there is no error and $\Rop^t(0)=0$. The stochastic operators $\Rop^t$ will be applied to a sequence of random vectors that will converge to zero, and hence the error will converge to zero as well, due to the relative variance property \eqref{eqrelva}. \algn{RandProx} 
 is therefore a \emph{variance-reduced} method \citep{han19,gor202,gow20a}: the random errors vanish along the iterations and the algorithm converges to an exact solution of the problem. 

To characterize how the error on the dual variable propagates to the primal variable after applying $K^*$, we also introduce the relative variance  $\oma \geq 0$ in the range of $K^*$ and the offset $\zeta \in [0,1]$ such that, for every $t\geq 0$,
\begin{equation}
\Exp{\sqnorm{K^*\big(\Rop^t(r^t)-r^t\big)}\;|\;\mathcal{F}_t} \leq \oma \sqnorm{r^t} - \zeta \sqnorm{K^*r^t}.\label{eqoma}
\end{equation}
It is easy to see that \eqref{eqoma} holds with $\oma=\|K\|^2\omega$ and $\zeta=0$, so this is the default choice without particular knowledge on $K^*$. But in some situations, e.g.\ sampling like in Section~\ref{secsamp}, a much smaller value of $\oma$ and a positive value of $\zeta$ can be derived.

\subsection{Description of the algorithm}

Let us now describe how the PDDY and \algn{RandProx} algorithms work. An iteration consists  in 3 steps: 
\begin{enumerate}
\item Given $x^{t}$ and $u^{t}$, the updated value of the primal variable is \emph{predicted} to be $\hat{x}^{t}$.
\item  The points $\hat{x}^{t}$ and $u^{t}$ are used to update the dual variable to its new value $u^{t+1}$.
\item The primal variable is \emph{corrected} from $\hat{x}^{t}$ to $x^{t+1}$, by back-propagating the difference $u^{t+1}-u^t$ using $K^*$.
\end{enumerate}
In \algn{RandProx}, randomization takes place in Step 2. On average, this decreases the progress from $u^{t}$ to $u^{t+1}$, and in turn from $\hat{x}^{t}$ to $x^{t+1}$ in Step 3, but the progress from $x^{t}$ to $\hat{x}^{t}$, due to the unaltered proximal gradient descent step in Step 1, is kept.  
Therefore, randomization can be used to balance the progress speed on the primal and dual variables, depending on the relative computational complexity of the gradient and proximity operators. The random errors are kept under control and convergence is ensured using \emph{underrelaxation}: let us define, for every $t\geq 0$,
\begin{equation}
\hat{u}^{t+1}\coloneqq \mathrm{prox}_{\tau \hc} \big(u^t+\tau K \hat{x}^{t}\big).\label{eqkolgkr}
\end{equation}
The PDDY algorithm updates the dual variable by setting $u^{t+1}\coloneqq \hat{u}^{t+1}$. In \algn{RandProx}, let us define $$\tilde{u}^{t+1}\coloneqq u^t +\Rop^t\big(\hat{u}^{t+1}-u^t\big)=\hat{u}^{t+1}+e^t$$ for some zero-mean random error $e^t$, keeping in mind that $\tilde{u}^{t+1}$ is typically cheaper to compute than $\hat{u}^{t+1}$. Then underrelaxation is applied: we set 
\begin{equation}
u^{t+1}\coloneqq \rho \tilde{u}^{t+1}+(1-\rho)u^t
\end{equation}
for some relaxation parameter $\rho \in (0,1]$; we use $\rho = \frac{1}{1+\omega}$ in the  algorithm. 
That is, the update of the dual variable consists in
a convex combination of the old estimate $u^t$ and the new, better in expectation but noisy, estimate $\tilde{u}^{t+1}$.  Noise is mitigated by underrelaxation, because the error $e^t$ is multiplied by $\rho$, so that its variance is multiplied by $\rho^2$. So, even if $\omega$ is arbitrarily large, $\omega\rho^2$ is kept small. Underrelaxation slows down the progress on the dual variable of the algorithm towards the solution, but if the iterations become faster, this is beneficial overall.

\section{Convergence analysis of \algn{RandProx}}

Our most general result, whose proof is in the Appendix, is the following:\medskip

\noindent\textbf{Theorem 1.}\ \  \emph{Suppose that $\mu_{\ff}>0$ or $\mu_{\g}>0$, and that $\mu_{\hc}>0$.
In \algn{RandProx}, suppose that $0<\gamma< \frac{2}{L_{\ff}}$, $\tau>0$, and $\gamma \tau \big((1-\zeta)\|K\|^2+\oma\big)\leq 1$, where $\oma$ and $\zeta$ are defined in \eqref{eqoma}.\footnote{The condition $\gamma < \frac{2}{L_{\ff}}$ is given for simplicity.  Larger values of $\gamma$ can be used when $\mu_{\g}>0$, as long as $c<1$ in \eqref{eqrate1}.}
For every $t\geq 0$, define the Lyapunov function
\begin{equation}
\Psi^{t}\coloneqq  \frac{1}{\gamma}\sqnorm{x^{t}-x^\star}+(1+\omega)\left(\frac{1}{\tau}+2\mu_{\hc}\right)\sqnorm{u^{t}-u^\star},\label{eqlya1p}
\end{equation}
where $x^\star$ and $u^\star$ are the unique solutions to \eqref{eqpb0} and \eqref{eqpbd}, respectively.
Then  \algn{RandProx} converges linearly:  for every $t\geq 0$,
\begin{align}
\Exp{\Psi^{t}}&\leq c^t \Psi^0,
\end{align}
where 
\begin{equation}
c\coloneqq  \max\left(\frac{(1-\gamma\mu_{\ff})^2}{1+\gamma\mu_{\g}},\frac{(\gamma L_{\ff}-1)^2}{1+\gamma\mu_{\g}},1-\frac{2\tau \mu_{\hc}}{(1+\omega)(1+2\tau \mu_{\hc})}\right)<1.\label{eqrate1}
\end{equation}
Also, $(x^t)_{t\in\mathbb{N}}$ and $(\hat{x}^t)_{t\in\mathbb{N}}$ both converge  to $x^\star$ and $(u^t)_{t\in\mathbb{N}}$ converges to $u^\star$, almost surely.
}\medskip

In Theorem 1, if $\gamma\leq \frac{2}{L_{\ff}+\mu_{\ff}}$, we have $\max(1-\gamma\mu_{\ff},\gamma L_{\ff}-1)^2=(1-\gamma\mu_{\ff})^2\leq 1-\gamma\mu_{\ff}$, so that in that case the rate $c$ in \eqref{eqrate1} satisfies
\begin{equation*}
c\leq 1-\min\left(\frac{\gamma(\mu_{\ff}+\mu_{\g})}{1+\gamma\mu_{\g}},\frac{2\tau \mu_{\hc} }{(1+\omega)(1+2\tau \mu_{\hc})}\right)<1.
\end{equation*}

\noindent \textbf{Remark 1} (choice of $\tau$)\ \ 
Given $\gamma$, the rate $c$ in \eqref{eqrate1} is smallest if $\tau$ is largest. So, there seems to be no reason to take $\tau\gamma \big((1-\zeta)\|K\|^2+\oma\big)<1$, and $\tau\gamma \big((1-\zeta)\|K\|^2+\oma\big)=1$ should be the best choice in most cases. Thus, one can set 
$\tau=\frac{1}{\gamma ((1-\zeta)\|K\|^2+\oma)}$
 and keep $\gamma$ as the only parameter to tune in \algn{RandProx}. \medskip

\begin{table}[t]
\centering \footnotesize
\caption{The different particular cases of the problem \eqref{eqpb0} for which we derive an instance of \algn{RandProx}, with the number of the theorem where its linear convergence is stated, and  the corresponding condition on $\h$ and $K$. $\lambda$ is a shorthand notation for $\lambda_{\min}(KK^*)$ and $\imath_{\{b\}}:x\mapsto (0$ if $x=b$, $+\infty$ otherwise$)$.}
\label{tab1}
\begin{tabular}{cccccccc}
\hline
$\ff$&$\g$&$\h$&$K$&\makecell{Deterministic\\algorithm}& \makecell{Randomized\\algorithm}&Theorem&\makecell{Condition ensuring\\linear convergence}\\
\hline
any&any&any&any&PDDY&\algn{RandProx}&1&$\mu_{\hc}>0$\\
any&0&any&any&PAPC&\algn{RandProx}&2&$\!\mu_{\hc}\!>\!0$ or $\lambda\!>\!0$\\
any&0&any&$\mathrm{Id}$&forward-backward (FB)&\algn{RandProx-FB}&3&---\\
any&0&$\imath_{\{b\}}$&any&PAPC&\algn{RandProx-LC}&4&---\\
0&any&any&any&Chambolle--Pock (CP)&\algn{RandProx-CP}&7&$\mu_{\hc}>0$\\
0&any&any&$\mathrm{Id}$&ADMM&\algn{RandProx-ADMM}&8&$\mu_{\hc}>0$\\
any&any&any&$\mathrm{Id}$&Davis--Yin (DY)&\algn{RandProx-DY}&9&$\mu_{\hc}>0$\\
\hline
\end{tabular}
\end{table}

In the rest of this section, we discuss some particular cases of \eqref{eqpb0}, for which we derive stronger convergence guarantees than in Theorem 1 for \algn{RandProx}. Other particular cases are studied in the Appendix; for instance, an instance of \algn{RandProx}, called \algn{RandProx-ADMM}, is a randomized version of the popular ADMM \citep{boy11}. The different particular cases are summarized in Table~\ref{tab1}.

\subsection{Particular case $\g=0$}\label{secpapc}

In this section, we assume that $\g=0$. Then the PDDY algorithm becomes an algorithm proposed for least-squares problems~\citep{lor11} and rediscovered independently as the PDFP2O algorithm~\citep{che13} and as the Proximal Alternating Predictor-Corrector (PAPC) algorithm~\citep{dro15}; let us call it the PAPC algorithm. It has been shown to have a primal--dual forward--backward structure \citep{com14}. Thus, when $\g=0$, \algn{RandProx} is a randomized version of the PAPC algorithm.

We note that $\ffc$ is strongly convex, which is not the case of $(\ff +\g)^*$ in general. Let us define $\lambda_{\min}(KK^*)$ as the smallest eigenvalue of $KK^*$.
 $\lambda_{\min}(KK^*)>0$ if and only if $\mathrm{ker}(K^*)=\{0\}$, where $\mathrm{ker}$ denotes the kernel. 
If $\lambda_{\min}(KK^*)>0$, $\ffc(-K^*\cdot)$ is strongly convex. Thus,  when $\g=0$, $\lambda_{\min}(KK^*)>0$ and $\mu_{\hc}>0$ are two sufficient conditions for the dual problem \eqref{eqpbd} to be strongly convex. We indeed get linear convergence of \algn{RandProx} in that case:\medskip

\noindent\textbf{Theorem 2.}\ \  \emph{Suppose that $\g=0$, $\mu_{\ff}>0$, and that $\lambda_{\min}(KK^*)>0$ or $\mu_{\hc}>0$.
In \algn{RandProx}, suppose that $0<\gamma < \frac{2}{L_{\ff}}$, $\tau>0$ and $\gamma \tau \big((1-\zeta)\|K\|^2+\oma\big)\leq 1$. 
Then \algn{RandProx} converges linearly:  for every $t\geq 0$,
$\Exp{\Psi^{t}}\leq c^t \Psi^0$,
where the Lyapunov function $\Psi^t$ is defined in \eqref{eqlya1p}, and 
\begin{equation}
c\coloneqq  \max\left((1-\gamma\mu_{\ff})^2,(\gamma L_{\ff}-1)^2,1-\frac{2\tau \mu_{\hc} + \gamma\tau\lambda_{\min}(KK^*)}{(1+\omega)(1+2\tau \mu_{\hc})}\right)<1.\label{eqrate2}
\end{equation}
Also, $(x^t)_{t\in\mathbb{N}}$ and $(\hat{x}^t)_{t\in\mathbb{N}}$ both converge  to $x^\star$ and $(u^t)_{t\in\mathbb{N}}$ converges to $u^\star$, almost surely.
}\medskip

When $\Rop^t=\mathrm{Id}$ and $\omega=\oma=0$, \algn{RandProx} reverts to the PAPC algorithm. Even in this particular case, Theorem 2 proves linear convergence of the PAPC algorithm and is new.  In \citet[Theorem 3.7]{che13}, the authors proved linear convergence of an underrelaxed version of the algorithm; underrelaxation slows down convergence.  In \citet{luk18}, Theorem 3.1 is wrong, since it is based on the false assumption that if $\lambda_{\min}(K_i K_i^*)>0$ for linear operators $K_i$, $i=1,\ldots,p$, then $\lambda_{\min}(K K^*)>0$, with $K:x\mapsto (K_1 x,\ldots,K_p x)$. Their theorem remains valid when $p=1$, but their rate is complicated and worse than ours.\smallskip

\begin{figure*}[t]
\begin{minipage}{.55\textwidth}
\begin{algorithm}[H]
		\caption{\algn{RandProx-FB} [new]}
		\begin{algorithmic}
			\STATE  \textbf{input:} initial points $x^0\in\mathcal{X}$, $u^0\in\mathcal{X}$;
			\STATE stepsize $\gamma>0$; $\omega\geq 0$
			\FOR{$t=0, 1, \ldots$}	
			\STATE $\hat{x}^{t} \coloneqq  x^t -\gamma \nabla \ff(x^t) - \gamma u^t$
			\STATE  $d^t\coloneqq \Rop^t \big( \hat{x}^{t} - \mathrm{prox}_{\gamma(1+\omega) \h} (\hat{x}^{t}+\gamma(1+\omega)u^t)\big)$
			\STATE  $u^{t+1}\coloneqq u^t +\frac{1}{\gamma(1+\omega)^2}d^t$
			\STATE $x^{t+1} \coloneqq \hat{x}^{t}-\frac{1}{1+\omega} d^t$
			\ENDFOR
		\end{algorithmic}
	\end{algorithm}\end{minipage}\ \ \ \ \ \begin{minipage}{.43\textwidth}
	\begin{algorithm}[H]
		\caption{\algn{RandProx-LC} [new]}
		\begin{algorithmic}
			\STATE  \textbf{input:} initial points $x^0\in\mathcal{X}$, $u^0\in\mathcal{U}$; 
			\STATE stepsizes $\gamma>0$, $\tau>0$; $\omega\geq 0$
			\STATE $v^0\coloneqq K^* u^0$
			\FOR{$t=0, 1, \ldots$}
			\STATE $\hat{x}^{t} \coloneqq  x^t -\gamma \nabla \ff(x^t) - \gamma v^t$
			\STATE $u^{t+1}\coloneqq u^t +\frac{\tau}{1+\omega}\Rop^t(K \hat{x}^{t}-b)$
			 \STATE $v^{t+1}\coloneqq K^* u^{t+1}$
			\STATE $x^{t+1} \coloneqq \hat{x}^{t}-\gamma (1+\omega) (v^{t+1}-v^t)$
			\ENDFOR
		\end{algorithmic}
	\end{algorithm}\end{minipage}\end{figure*}

We now consider the even more particular case of $\g=0$ and $K=\mathrm{Id}$. Then the problems \eqref{eqpb0} and \eqref{eqpbd} consist in minimizing $\ff(x)+\h(x)$ and $\ffc(-u)+\hc(u)$, respectively. The dual problem is strongly convex and has a unique solution $u^\star = -\nabla \ff(x^\star)$, for any primal solution $x^\star$. By setting $\tau \coloneqq 1/\gamma$ in the PAPC algorithm, we obtain the classical proximal gradient, a.k.a.\ forward-backward (FB), algorithm, which iterates $x^{t+1} \coloneqq \mathrm{prox}_{\gamma \h} \big(x^t -\gamma \nabla \ff(x^t)\big)$. Thus, when randomness is introduced, we set $\oma\coloneqq\omega$, $\zeta\coloneqq 0$ and, according to Remark 1, $\tau\coloneqq \frac{1}{\gamma(1+\omega)}$ in \algn{RandProx}. By noting that, for every $a>0$, the abstract operators $\Rop^t$ and  $a\Rop^t\big(\frac{1}{a} \cdot\big)$ have the same properties, we can put the constant $\gamma (1+\omega)$ outside $\Rop^t$ to simplify the algorithm, and rewrite \algn{RandProx} as \algn{RandProx-FB}, shown above. 
As a corollary of Theorem 2, we have:\medskip

\noindent\textbf{Theorem 3.}\ \  \emph{Suppose that $\mu_{\ff}>0$. In \algn{RandProx-FB}, 
suppose that $0<\gamma < \frac{2}{L_{\ff}}$. 
For every $t\geq 0$, define the Lyapunov function
\begin{equation}
\Psi^{t}\coloneqq  \frac{1}{\gamma}\sqnorm{x^{t}-x^\star}+(1+\omega)\big(\gamma(1+\omega)+2\mu_{\hc}\big)\sqnorm{u^{t}-u^\star},\label{eqlya1j}
\end{equation}
where $x^\star$ is the unique minimizer of $\ff+\h$ and $u^\star = -\nabla \ff(x^\star)$ is the unique minimizer of $\ffc(-\cdot)+\hc$.  Then 
\algn{RandProx-FB}  
converges linearly:  for every $t\geq 0$, 
$$\Exp{\Psi^{t}}\leq c^t \Psi^0,$$ 
where 
\begin{equation}
c\coloneqq  \max\left((1-\gamma\mu_{\ff})^2,(\gamma L_{\ff}-1)^2,1-\frac{1+\frac{2}{\gamma}\mu_{\hc}}{(1+\omega)\big(1+\omega+\frac{2}{\gamma} \mu_{\hc}\big)}\right)<1.\label{eqrate2j}
\end{equation}
Also, $(x^t)_{t\in\mathbb{N}}$ and $(\hat{x}^t)_{t\in\mathbb{N}}$ both converge  to $x^\star$ and $(u^t)_{t\in\mathbb{N}}$ converges to $u^\star$, almost surely.
}\medskip

It is important to note that it is not necessary to have $\mu_{\hc}>0$ in Theorem 3. If we ignore the properties of $\hc$, the third factor in \eqref{eqrate2j} can be replaced by its upper bound 
$
1-\frac{1}{(1+\omega)^2}.
$

\subsection{Linearly constrained smooth minimization}\label{seclcm}

Let $b\in\mathrm{ran}(K)$. In this section, we consider the linearly constrained (LC) minimization problem
\begin{equation}
\mathrm{Find} \ x^\star \in \argmin_{x\in\mathcal{X}} \, \ff(x) \quad \mbox{s.t.} \quad Kx=b,\label{eqlck1}
 \end{equation}
which is a particular case of \eqref{eqpb0} with $\g=0$ and 
$\h:u\in\mathcal{U}\mapsto(0$ if $u=b$, $+\infty$ otherwise$)$. 
We have $\hc:u\in\mathcal{U}\mapsto\langle u,b\rangle $ and $\mathrm{prox}_{\tau \hc}:u\in\mathcal{U}\mapsto u-\tau b$. 
The dual problem to \eqref{eqlck1} is
\begin{equation}
\mathrm{Find} \ u^\star \in \argmin_{u\in\mathcal{U}} \, \Big(\ffc(-K^*u)+\langle u,b\rangle\Big).\label{eqpbd2}
 \end{equation}
 We denote by $u^{\star}_0$ the unique solution to \eqref{eqpbd2} in $\mathrm{ran}(K)$. Then the set of solutions of \eqref{eqpbd2} is the affine subspace $u^{\star}_0+\mathrm{ker}(K^*)$. Thus, the dual problem is not strongly convex, unless $\mathrm{ker}(K^*)=\{0\}$. Yet, we will see that strong convexity of $\ff$ is sufficient to have linear convergence of \algn{RandProx},  without any condition on $K$.
 
 We rewrite  \algn{RandProx} in this setting as \algn{RandProx-LC}, shown above. 
We observe that $u^t$ does not appear in the argument of $\Rop^t$ any more, so that the iteration can be rewritten with the variable $v^t=K^*u^t$, and $u^t$ can be removed if we are not interested in estimating a dual solution. In any case, we denote by $P_{\mathrm{ran}(K)}$ the orthogonal projector onto $\mathrm{ran}(K)$ and by $\lambda_{\min}^+(KK^*)>0$  the smallest \emph{nonzero} eigenvalue of $KK^*$. Then:\medskip

\noindent\textbf{Theorem 4.}\ \  \emph{In the setup \eqref{eqlck1}--\eqref{eqpbd2}, suppose that $\mu_{\ff}>0$. In \algn{RandProx-LC}, suppose that $0<\gamma < \frac{2}{L_{\ff}}$, $\tau>0$ and $\gamma \tau \big((1-\zeta)\|K\|^2+\oma\big)\leq 1$. 
Define the Lyapunov function, for every $t\geq 0$,
\begin{equation}
\Psi^{t}\coloneqq  \frac{1}{\gamma}\sqnorm{x^{t}-x^\star}+\frac{1+\omega}{\tau}\sqnorm{u^{t}_0-u_0^\star},
\end{equation}
where $u^t_0\coloneqq P_{\mathrm{ran}(K)}(u^t)$ is also the unique element in $\mathrm{ran}(K)$ such that $v^t=K^*u^t_0$, $x^\star$ is the unique solution of \eqref{eqlck1} and $u_0^\star$ is the unique solution in $\mathrm{ran}(K)$ of \eqref{eqpbd2}. Then  \algn{RandProx-LC} converges linearly:  for every $t\geq 0$, 
$$\Exp{\Psi^{t}}\leq c^t \Psi^0,$$
where 
\begin{equation}
c\coloneqq \max\left((1-\gamma\mu_{\ff})^2,(\gamma L_{\ff}-1)^2,1-\frac{\gamma\tau\lambda_{\min}^+(KK^*)}{1+\omega}\right)<1.\label{eqratelc}
\end{equation}
Also, $(x^t)_{t\in\mathbb{N}}$ and $(\hat{x}^t)_{t\in\mathbb{N}}$ both converge  to $x^\star$ and $(u_0^t)_{t\in\mathbb{N}}$ converges to $u_0^\star$, almost surely.
}
\medskip

Theorem 4 is new even for the PAPC algorithm when $\omega=0$: its linear convergence under the stronger condition $\gamma \tau \|K\|^2< 1$ has been shown in \citet[Theorem 6.2]{sal20}, but our rate in \eqref{eqratelc} is better.

We further discuss \algn{RandProx-LC}, which can be used for decentralized optimization, in the Appendix. Another example of application is  when $\mathcal{X}=\mathbb{R}^d$, for some $d\geq 1$, and $K$ is a matrix; one can solve \eqref{eqlck1} by activating one row of $K$ chosen uniformly at random at every iteration.

\section{Convergence in the merely convex case}

In all theorems, strong convexity of $\ff$ or $\g$ is assumed; that is, $\mu_{\ff}>0$ or $\mu_\g>0$. In this section, we remove this hypothesis, so that the primal problem is not necessarily strongly convex any more. But $\nabla \ff (x^\star)$ is the same for every solution $x^\star$ of \eqref{eqpb0}, and we denote by $\nabla \ff (x^\star)$ this element.

We define the Bregman divergence of $\ff$ at points $(x,x')\in\mathcal{X}^2$ as $$D_{\ff}(x,x')\coloneqq \ff(x)-\ff(x')-\langle \nabla \ff (x'),x-x'\rangle \geq 0.$$ For every $t\geq 0$, $D_{\ff}(x^t,x^\star)$  is the same for every solution $x^\star$ of \eqref{eqpb0}, and we denote by $D_{\ff}(x^t,x^\star)$ this element. $D_{\ff}(x^t,x^\star)$ can be viewed as a generalization of the objective gap $\ff(x^t) -\ff(x^\star)$ to the case when $\nabla \ff(x^\star)\neq 0$. $D_{\ff}(x^t,x^\star)$ is a loose kind of distance between $x^t$ and the solution set, but under some additional assumptions on $\ff$,  for instance strict convexity, $D_{\ff}(x^t,x^\star)\rightarrow 0$ implies that the distance from $x^t$ to the solution set tends to zero. Also, $D_{\ff}(x^t,x^\star)\geq \frac{1}{2L_\ff} \| \nabla \ff(x^t)-\nabla \ff(x^\star)\|^2$, so that $D_{\ff}(x^t,x^\star)\rightarrow 0$ implies that $\big(\nabla \ff(x^t)\big)_{t\in\mathbb{N}}$ converges to  $\nabla \ff(x^\star)$.
\medskip

\noindent\textbf{Theorem 11.}\ \  {\itshape In \algn{RandProx}, suppose that $0<\gamma< \frac{2}{L_{\ff}}$, $\tau>0$, and $\gamma \tau \big((1-\zeta)\|K\|^2+\oma\big)\leq 1$. Then $D_{\ff}(x^t,x^\star)\rightarrow 0$, almost surely and in quadratic mean. Moreover, for every $t\geq 0$, we define $\bar{x}^t \coloneqq \frac{1}{t+1} \sum_{i=0}^t x^i$. Then, for every $t\geq 0$,
\begin{equation}
\Exp{D_\ff(\bar{x}^t,x^\star)} \leq \frac{\Psi^0}{(2-\gamma L_\ff)(t+1)}=\mathcal{O}(1/t).
\end{equation}

If, in addition, $\mu_{\hc}>0$, there is a unique dual solution $u^\star$ to \eqref{eqpbd} and $(u^t)_{t\in\mathbb{N}}$ converges to $u^\star$, in quadratic mean.
}\medskip

We can derive counterparts of the other theorems in the same way. These theorems apply to all algorithms presented in the paper. For instance, Theorem 11 applies to Scaffnew~\citep{mis22}, a particular case of \algn{RandProx-FL} seen in Section~\ref{secfl}, 
and provides for it the first convergence results in the non-strongly convex case.

\subsubsection*{Acknowledgments}
The work of P. Richt\'{a}rik was partially supported by the KAUST Baseline Research Fund Scheme and by the SDAIA-KAUST Center of Excellence in Data Science and Artificial Intelligence.

\bibliographystyle{iclr2023_conference}
\bibliography{IEEEabrv,biblio}

\clearpage
\appendix

{\huge Appendix}

\section{Examples}

\subsection{Skipping the proximity operator}\label{secskip}

\begin{figure*}[t]
\begin{minipage}{.44\textwidth}
\begin{algorithm}[H]
		\caption{\algn{RandProx-Skip} [new]}
		\begin{algorithmic}
			\STATE  \textbf{input:} initial points $x^0\in\mathcal{X}$, $u^0\in\mathcal{U}$; \\
			stepsizes $\gamma>0$, $\tau>0$; $p\in (0,1]$
			\STATE $v^0\coloneqq K^* u^0$
			\FOR{$t=0, 1, \ldots$}
			\STATE $\hat{x}^{t} \coloneqq  \mathrm{prox}_{\gamma \g}\big(x^t -\gamma \nabla \ff(x^t) - \gamma v^t\big)$
			\STATE Flip a coin $\theta^t=(1$ with probability $p$, $0$ else$)$
			\IF{$\theta^t=1$}
			\STATE $u^{t+1}\coloneqq \mathrm{prox}_{\tau \hc} (u^t+\tau K \hat{x}^{t})$
			\STATE $v^{t+1}\coloneqq K^* u^{t+1}$
			\STATE  $x^{t+1} \coloneqq \hat{x}^{t}-\frac{\gamma}{p}(v^{t+1}-v^t)$
			\ELSE
			\STATE $u^{t+1}\coloneqq u^t$, $v^{t+1}\coloneqq v^t$, $x^{t+1} \coloneqq \hat{x}^{t}$
			\ENDIF
			\ENDFOR
		\end{algorithmic}
	\end{algorithm}\end{minipage}\ \ \ \ \ \begin{minipage}{.52\textwidth}
	\begin{algorithm}[H]
		\caption{\algn{RandProx-Minibatch} [new]}
		\begin{algorithmic}
			\STATE  \textbf{input:} initial points $x^0\in\mathcal{X}$, $(u_i^0)_{i=1}^n\in\mathcal{X}^n$; 
			\STATE stepsize $\gamma>0$; $k\in \{1,\dots,n\}$
			\STATE $v^0\coloneqq \sum_{i=1}^n u_i^0$
			\FOR{$t=0, 1, \ldots$}
			\STATE $\hat{x}^{t} \coloneqq  \mathrm{prox}_{\gamma \g}\big(x^t -\gamma \nabla \ff(x^t) - \gamma v^t\big)$
			\STATE pick $\Omega^t\subset \{1,\ldots,n\}$ of size $k$ unif.\ at random
			\FOR{$i\in\Omega^t$}
			\STATE $u_{i}^{t+1}\coloneqq \mathrm{prox}_{\frac{1}{\gamma n} \hic} (u_{i}^t+\frac{1}{\gamma n}  \hat{x}^{t})$
			\ENDFOR
			\FOR{$i\in\{1, \ldots,n\}\backslash \Omega^t$}
			\STATE $u_{i}^{t+1}\coloneqq u_{i}^t $
			\ENDFOR
			\STATE $v^{t+1}\coloneqq \sum_{i=1}^n u_i^{t+1}$
			\STATE  $x^{t+1} \coloneqq \hat{x}^{t}- \frac{\gamma n}{k} (v^{t+1}-v^t)$
			\ENDFOR
		\end{algorithmic}
	\end{algorithm}\end{minipage}\end{figure*}

In this section, we consider the case of Bernoulli operators $\Rop^t$ defined in \eqref{eqber}, which compute and return their argument only with probability $p>0$. \algn{RandProx} becomes \algn{RandProx-Skip}, shown above. Then $\omega = \frac{1}{p}-1$, $\oma=\|K\|^2\omega$, and $\zeta=0$.

If $\g=0$, \algn{RandProx-Skip} reverts to the SplitSkip algorithm proposed recently \citep{mis22}. Our Theorems 1 and 4 recover the same rate as given for SplitSkip in \citet[Theorem D.1]{mis22}, if smoothness of $\h$ is ignored. If in addition $K=\mathrm{Id}$ and $\tau=\frac{1}{\gamma(1+\omega)}=\frac{p}{\gamma}$, \algn{RandProx-Skip} reverts to ProxSkip, a particular case of SplitSkip \citep{mis22}. Our Theorem 3 applies to this case and allows us to exploit the possible smoothness of $\h$ in \algn{RandProx-Skip}=\,ProxSkip, which is not the case of the results of \citep{mis22}. As a practical application of our new results, let us consider \emph{personalized federated learning (FL)} \citep{han20}: given 
a client-server architecture with a master and $n\geq 1$ users, each with local cost function $\ffi$, $i=1,\ldots,n$, the goal is to
\begin{equation}
 \minimize_{(x_i)_{i=1}^n\in (\mathbb{R}^d)^n}\;  \sum_{i=1}^n\ffi(x_i)+\frac{\lambda}{2} \sum_{i=1}^n \|x_i- \bar{x}\|^2,
\end{equation}
where $\bar{x}\coloneqq \frac{1}{n} \sum_{i=1}^n x_i$. Each $\ffi$ is supposed $L_{\ff}$-smooth and $\mu_{\ff}$-strongly convex. 
We set $\mathcal{X}\coloneqq(\mathbb{R}^d)^n$,  $\ff:x=(x_i)_{i=1}^n\mapsto \sum_{i=1}^n\ffi(x_i)$, $\h:x\mapsto \frac{\lambda}{2} \sum_{i=1}^n \|x_i- \bar{x}\|^2$. $\ff$ is $L_{\ff}$-smooth and $\mu_{\ff}$-strongly convex, $\h$  is $\lambda$-smooth, so  that $\mu_{\hc}=\frac{1}{\lambda}$. Thus,  with $\gamma= \frac{1}{L_{\ff}}$, we have in \eqref{eqrate2j}:
\begin{equation*}
c\leq   1-\min\left(\frac{\mu_{\ff}}{L_{\ff}},\frac{1+\frac{2L_{\ff}}{\lambda}}{\frac{1}{p}\big(\frac{1}{p}+\frac{2L_{\ff}}{\lambda}\big)}\right)<1.
\end{equation*}
Hence, with $p=\frac{\sqrt{\mu_{\ff}\min(L_{\ff},\lambda)}}{L_{\ff}}=\sqrt{\frac{\mu_{\ff}}{L_{\ff}}}\min\Big(\sqrt{\frac{\lambda}{L_{\ff}}},1\Big)$, the communication complexity in terms of the expected number of communication rounds to reach $\epsilon$-accuracy is $\mathcal{O}\left(\left(\sqrt{\frac{\min(L_{\ff},\lambda)}{\mu_{\ff}}}+1\right)\log\frac{1}{\epsilon}\right)$, which, up to the `+1' log factor, is optimal \citep{han20}. This shows that in personalized FL with $\lambda<L_{\ff}$, the complexity can be decreased in comparison with non-personalized FL, which corresponds to $\lambda=+\infty$. This is achieved by properly setting $p$ in ProxSkip, according to our new theory, which exploits the smoothness of  $\h$.

\subsection{Sampling among several functions}\label{secsamp}

We first remark that we can extend Problem \eqref{eqpb0} with the term $\h(Kx)$ replaced by the sum $\sum_{i=1}^n \hi(K_i x)$ of $n\geq 2$ proper closed convex functions $\hi$ composed with linear operators $K_i:\mathcal{X}\rightarrow \mathcal{U}_i$, for some real Hilbert spaces $\mathcal{U}_i$, by using the classical product-space trick: by defining $\mathcal{U}\coloneqq \mathcal{U}_1 \times \cdots\mathcal{U}_n$, $\h:u=(u_i)_{i=1}^n \in \mathcal{U} \mapsto \sum_{i=1}^n \hi(u_i)$, $K: x\in\mathcal{X}\mapsto   (K_ix)_{i=1}^n \in\mathcal{U}$, we have $\h(Kx)=\sum_{i=1}^n \hi(K_i x)$.    
In particular, by setting $K_i\coloneqq \mathrm{Id}$  and $\mathcal{U}_i\coloneqq \mathcal{X}$, 
we consider in this section the problem:
\begin{equation}
\mathrm{Find} \ x^\star \in \argmin_{x\in\mathcal{X}}  \left( \ff(x) + \g(x)+\sum_{i=1}^n \hi(x)\right).\label{eqpbm}
 \end{equation}
 We have $\hc:(u_i)_{i=1}^n\in \mathcal{X}^n\mapsto \sum_{i=1}^n \hic(u_i)$ and we suppose that every function $\hic$ is $\mu_{\hc}$-strongly convex, for some $\mu_{\hc}\geq 0$; then $\hc$ is $\mu_{\hc}$-strongly convex.
  Thus, the dual problem to \eqref{eqpbm} is
 \begin{equation}
\mathrm{Find} \ (u_i^\star)_{i=1}^n \in \argmin_{(u_i)_{i=1}^n\in\mathcal{X}^n} \, \left( (\ff+\g)^*\Big({\textstyle-\sum \limits_{i=1}^n u_i}\Big)+\sum_{i=1}^n \hic(u_i)\right).\label{eqpbmd}
 \end{equation}
 Since $K^*K = n\mathrm{Id}$, $\|K\|^2 = n$.  Now, we choose  $\Rop^t$ as the \texttt{rand}-$k$ sampling operator, for some  $k \in \{1,\ldots,n\}$: $\Rop^t$  multiplies $k$ elements out of the $n$ of its argument sequence, chosen uniformly at random, by $n/k$ and sets the other ones to zero. It is known \citep[Proposition 1]{con22m} that we can set 
 \begin{equation*}
 \omega \coloneqq \frac{n}{k}-1,\quad \oma \coloneqq \frac{n(n-k)}{k(n-1)},\quad \zeta\coloneqq \frac{n-k}{k(n-1)}.
 \end{equation*}
 Note that this value of $\oma$ is $n-1$ times smaller than the naive bound $\|K\|^2\omega = \frac{n(n-k)}{k}$. We have 
 $(1-\zeta)\|K\|^2 + \oma=n$. \algn{RandProx} in this setting, with $\tau\coloneqq \frac{1}{\gamma n}$, becomes \algn{RandProx-Minibatch}, shown above, and Theorem 1 yields:\medskip

\noindent\textbf{Theorem 5.}\ \  \emph{Suppose that $\mu_{\ff}>0$ or $\mu_{\g}>0$, and that $\mu_{\hc}>0$.
In \algn{RandProx-Minibatch}, suppose that $0<\gamma< \frac{2}{L_{\ff}}$.  
Define the Lyapunov function, for every $t\geq 0$,
\begin{equation}
\Psi^{t}\coloneqq  \frac{1}{\gamma}\sqnorm{x^{t}-x^\star}+\frac{n}{k}\left(\gamma n+2\mu_{\hc}\right)\sum_{i=1}^n\sqnorm{u_i^{t}-u_i^\star},\label{eqlya4}
\end{equation}
where $x^\star$ and $(u_i^\star)_{i=1}^n$ are the unique solutions to \eqref{eqpbm} and \eqref{eqpbmd}, respectively.
Then \algn{RandProx-Minibatch} converges linearly:  for every $t\geq 0$, 
$\Exp{\Psi^{t}}\leq c^t \Psi^0$, 
where 
\begin{equation}
c\coloneqq  \max\left(\frac{(1-\gamma\mu_{\ff})^2}{1+\gamma\mu_{\g}},\frac{(\gamma L_{\ff}-1)^2}{1+\gamma\mu_{\g}},1-\frac{2k \mu_{\hc} }{n(\gamma n+2 \mu_{\hc})}\right).
\end{equation}
Also, $(x^t)_{t\in\mathbb{N}}$ and $(\hat{x}^t)_{t\in\mathbb{N}}$ both converge  to $x^\star$ and $(u_i^t)_{t\in\mathbb{N}}$ converges to $u_i^\star$, $\forall i$, 
almost surely.
}

 \algn{RandProx-Minibatch} with $k=1$ becomes the Stochastic Decoupling Method (SDM) proposed in \citet{mis19}, where strong convexity of $\g$ is not exploited, but similar guarantees are derived as in Theorem 5 if $\mu_{\g}=0$. 
 Linear convergence of SDM is also proved in \citet{mis19} in conditions related to ours in Theorems 2 and 4. 
 Thus, \algn{RandProx-Minibatch} extends SDM to larger minibatch size $k$ and exploits possible strong convexity of $\g$.

When $\ff=0$ and $\g=0$, SDM further simplifies to Point-SAGA \citep{def16}. In that case, our results do not apply directly, since there is no strong convexity in $\ff$ and $\g$ any more, but when minimizing the average of functions  $\hi$, with each function supposed to be $L$-smooth and $\mu$-strongly convex, for some $L\geq \mu >0$, we can transfer the strong convexity to $\g$ by subtracting $\frac{\mu}{2} \|\cdot\|^2$ to each $\hi$ and setting $\g=\frac{\mu}{2} \|\cdot\|^2$. This does not change the problem and the algorithm but our Theorem 5 now applies, and with the right choice of $\gamma$, we recover the result in \citet{def16}, that the asymptotic complexity of Point-SAGA to reach $\epsilon$-accuracy is $\mathcal{O}\left(\Big(n+\sqrt{\frac{nL}{\mu}}\Big)\log\frac{1}{\epsilon}\right)$, which is conjectured to be optimal. 

 Thus, \algn{RandProx-Minibatch} extends Point-SAGA to larger minibatch size  and to the more general problem  \eqref{eqpbm} with nonzero $\ff$ or $\g$.
 
 When $n=1$, there is no randomness and SDM reverts to the DY algorithm discussed in Appendix~\ref{seccdy}.
 
 \begin{figure*}[t]
\begin{minipage}{.48\textwidth}
	\begin{algorithm}[H]
		\caption{
		SDM \\ \citep{mis19}}
		\begin{algorithmic}
			\STATE  \textbf{input:} initial points $x^0\in\mathcal{X}$, $(u_i^0)_{i=1}^n\in\mathcal{X}^n$;
			\STATE stepsize $\gamma>0$
			\STATE $v^0\coloneqq \sum_{i=1}^n u_i^0$
			\FOR{$t=0, 1, \ldots$}
			\STATE $\hat{x}^{t} \coloneqq  \mathrm{prox}_{\gamma \g}\big(x^t -\gamma \nabla \ff(x^t) - \gamma v^t\big)$
			\STATE pick $i^t \in \{1,\ldots,n\}$  uniformly at random
			\STATE $x^{t+1}\coloneqq \mathrm{prox}_{\gamma n \hi} (\gamma n u_{i^t}^t+\hat{x}^{t})$
			\STATE $u_{i^t}^{t+1}\coloneqq u_{i^t}^t+\frac{1}{\gamma n} (\hat{x}^{t}-x^{t+1})$
			\STATE for every $i\in\{1, \ldots,n\}\backslash\{i^t\}$, $u_{i}^{t+1}\coloneqq u_{i}^t $
			\STATE $v^{t+1}\coloneqq \sum_{i=1}^n u_i^{t+1}$ // $=v^t + u_{i^t}^{t+1}-u_{i^t}^{t}$
			\ENDFOR
		\end{algorithmic}
\end{algorithm}\end{minipage}\ \ \ \ \ 
	\begin{minipage}{.48\textwidth}		
\begin{algorithm}[H]
		\caption{
		Point-SAGA \\ \citep{def16}}
		\begin{algorithmic}
			\STATE  \textbf{input:} initial points $x^0\in\mathcal{X}$, $(u_i^0)_{i=1}^n\in\mathcal{X}^n$;
			\STATE stepsize $\gamma>0$
			\STATE $v^0\coloneqq \sum_{i=1}^n u_i^0$
			\FOR{$t=0, 1, \ldots$}
			\STATE $\hat{x}^{t} \coloneqq  x^t - \gamma v^t$
			\STATE pick $i^t \in \{1,\ldots,n\}$  uniformly at random
			\STATE $x^{t+1}\coloneqq \mathrm{prox}_{\gamma n \hi} (\gamma n u_{i^t}^t+\hat{x}^{t})$
			\STATE $u_{i^t}^{t+1}\coloneqq u_{i^t}^t+\frac{1}{\gamma n} (\hat{x}^{t}-x^{t+1})$
			\STATE for every $i\in\{1, \ldots,n\}\backslash\{i^t\}$, $u_{i}^{t+1}\coloneqq u_{i}^t $
			\STATE $v^{t+1}\coloneqq \sum_{i=1}^n u_i^{t+1}$ // $=v^t + u_{i^t}^{t+1}-u_{i^t}^{t}$
			\ENDFOR
		\end{algorithmic}
		\end{algorithm}\end{minipage}\end{figure*}

\subsection{Distributed and federated learning with compression}\label{secfl}

\begin{figure*}[t]
\centering
\begin{minipage}{.55\textwidth}
\begin{algorithm}[H]
		\caption{\algn{RandProx-FL} [new]}
		\begin{algorithmic}
			\STATE  \textbf{input:} initial estimates $(x_i^0)_{i=1}^n \in\mathcal{X}^n$, $(u_i^0)_{i=1}^n\in\mathcal{X}^n$ such that  $\sum_{i=1}^n u_i^0 = 0$;  		
			stepsize $\gamma>0$; $\omega\geq 0$
			\FOR{$t=0, 1, \ldots$}	
			\FOR{$i=1, \ldots,n$ at nodes in parallel}
			\STATE $\hat{x}_i^{t} \coloneqq  x_i^t -\gamma \nabla \ffi(x_i^t) - \gamma u_i^t$
			\STATE $a_i^t \coloneqq   \Rop^t (\hat{x}_i^t )$
			\STATE // send compressed vector $a_i^t$ to master
			\ENDFOR
			\STATE $a^{t}\coloneqq \frac{1}{n}\sum_{i=1}^n a_i^t$\ \  \ // aggregation at master
			\STATE // broadcast $a^{t}$ to all nodes
			\FOR{$i=1, \ldots,n$ at nodes in parallel}
			\STATE $d_i^t \coloneqq  a_i^t - a^t$
			\STATE $u_i^{t+1}\coloneqq u_i^t + \frac{1}{\gamma(1+\omega)^2}d_i^t$
			\STATE $x_i^{t+1} \coloneqq \hat{x}_i^{t} - \frac{1}{1+\omega}d_i^t$
			\ENDFOR
			\ENDFOR
\end{algorithmic}
\end{algorithm}\end{minipage}
\end{figure*}

We consider in this section distributed optimization within the 
client-server model, with a master node communicating back and forth with $n\geq 1$ parallel workers.  This is particularly relevant for federated learning (FL)~\citep{kon16,mcm17,kai19,li20}, 
where a potentially huge number of devices, with their owners' data stored on each of them, are involved in the collaborative process of  training a global machine learning model. The goal is to exploit the wealth of useful information lying in the heterogeneous data stored across the devices. 
	Communication between the devices and the distant server, which can be costly and slow, is the main bottleneck in this framework. So, it is of primary importance to devise novel algorithmic strategies, which are efficient in terms of computation and communication complexities. A natural and widely used idea is 
	 to make use of (lossy) {\em compression}, to reduce the size of the communicated message~\citep{ali17,wen17,wan18,GDCI,alb20,bas20,dut20,sat20,xu21}. Another popular idea is to make use of \emph{local steps}~\citep{mcm17,localGD,sti19,kha20local,mal20,woo20,kar20,gor21,mis22}; 
	 that is, communication with the server does not occur at every iteration but only every few iterations, for instance communication is triggered randomly with a small probability at every iteration. Between communication rounds,  the workers perform multiple local steps  independently, based on their local objectives. Our proposed algorithm \algn{RandProx-FL} unifies the two strategies,   in the sense that depending on the choice of the randomization process $\Rop^t$, we obtain  a method with local steps or with compression, or both. The combination of local training and compression has been further investigated in our follow-up work \citep{con22cs}, and partial participation in \citet{con23pp}.

Thus, we consider the problem
\begin{equation}
\mathrm{Find} \ x^\star \in \argmin_{x\in\mathbb{R}^d}  \left(\sum_{i=1}^n \ffi(x) \right),\label{eqcompress1}
 \end{equation}
 where $d\geq 1$ is the model dimension and $n\geq 1$ is the number of parallel workers, each having its own objective function $\ffi$. 
 Every function $\ffi:\mathbb{R}^d \rightarrow\mathbb{R}$ is $\mu$-strongly convex and $L$-smooth, for some $L\geq \mu>0$. We define $\kappa \coloneqq L/\mu$.
 
 Now, we can observe that \eqref{eqcompress1} can be recast as \eqref{eqpb0} with $K=\mathrm{Id}$, $\mathcal{U} = \mathcal{X}$, $\g=0$; that is, as the minimization of $\ff + \h$, as studied in Section~\ref{secpapc}, with 
 \begin{align}
 &\mathcal{X} = (\mathbb{R}^d)^n, \quad \ff:x=(x_i)_{i=1}^n \mapsto  \sum_{i=1}^n \ffi(x_i),\\
 & \h:x=(x_i)_{i=1}^n \mapsto (0\mbox{ if } x_1=\cdots=x_n,\ +\infty \mbox{ otherwise}).
 \end{align}
 We note that $\ff$ is $\mu$-strongly convex and $L$-smooth, and  $\mu_{\hc}=0$. 
Making these substitutions in \algn{RandProx-FB} yields \algn{RandProx-FL}, a distributed algorithm well suited for FL, shown above. In \algn{RandProx-FL}, randomization takes the form of \emph{linear} random unbiased operators $\Rop^t$ applied to the vectors sent to the server. 
Note that at every iteration, 
the same operator $\Rop^t$ is applied at every node; that is, its randomness is shared.  We can easily check that \algn{RandProx-FL} is an instance of \algn{RandProx-FB},  because of the linearity of the $\Rop^t$  and because the property $\sum_{i=1}^n u_i^t = 0$ is maintained at every iteration. Formally, $\Rop^t$ applied as a whole in \algn{RandProx-FB}  consists of $n$ copies of $\Rop^t$ applied individually at every node in \algn{RandProx-FL}, that is why we keep the same notation; in particular, the value of $\omega$  is the same in both interpretations. 

Interestingly, in \algn{RandProx-FL}, information about the functions $\ffi$ or their gradients is never communicated and is exploited completely locally. This is ideal in terms of privacy.

As an application of Theorem 3, we obtain:\medskip

\noindent\textbf{Theorem 10.}\ \  \emph{In \algn{RandProx-FL}, 
suppose that $0<\gamma < \frac{2}{L_{\ff}}$. 
Define the Lyapunov function, for every $t\geq 0$,
\begin{equation}
\Psi^{t}\coloneqq  \sum_{i=1}^n \left( \frac{1}{\gamma} \sqnorm{x_i^{t}-x^\star}+\gamma(1+\omega)^2\sqnorm{u_i^{t}-u_i^\star}\right),\label{eqlya1jpp}
\end{equation}
where $x^\star$ is the unique solution of \eqref{eqcompress1} and $u_i^\star \coloneqq -\nabla \ffi(x^\star)$.  Then 
\algn{RandProx-FL} 
converges linearly:  for every $t\geq 0$, 
$\Exp{\Psi^{t}}\leq c^t \Psi^0$, 
where 
\begin{equation}
c\coloneqq  \max\left((1-\gamma\mu_{\ff})^2,(\gamma L_{\ff}-1)^2,1-\frac{1}{(1+\omega)^2}\right)<1.\label{eqrate2jpp}
\end{equation}
Also, the $(x_i^t)_{t\in\mathbb{N}}$ and $(\hat{x}_i^t)_{t\in\mathbb{N}}$ all converge  to $x^\star$ and every $(u_i^t)_{t\in\mathbb{N}}$ converges to $u_i^\star$, almost surely.}\medskip

If $\Rop^t$ is the Bernoulli compressor we have seen before in \eqref{eqber} and in Section~\ref{secskip}, \algn{RandProx-FL} reverts to the Scaffnew algorithm proposed in \citet{mis22}, which communicates at every iteration with probability $p \in (0,1]$ and performs in average $1/p$ local steps between successive communication rounds.  We have $\omega = \frac{1}{p}-1$. The analysis of Scaffnew in Theorem 10 is the same as in \citet{mis22}.  
With $\gamma=\frac{1}{L}$, the iteration complexity of Scaffnew  is $\mathcal{O}\big((\kappa+\frac{1}{p^2})\log\frac{1}{\epsilon}\big)$, and since the algorithm communicates with probability $p$, its average communication complexity is $\mathcal{O}\big((p\kappa+\frac{1}{p})\log\frac{1}{\epsilon}\big)$. In particular, with $p=\frac{1}{\sqrt{\kappa}}$,  the average communication complexity of Scaffnew is $\mathcal{O}\big(\sqrt{\kappa} \log\frac{1}{\epsilon}\big)$.

We now propose a new algorithm with compressed communication: in  \algn{RandProx-FL} we choose, for every $t\geq 0$, $\Rop^t$ as the well-known  \texttt{rand}-$k$ compressor, for some $k\in \{1,\ldots,d\}$: $\Rop^t$  multiplies $k$ coordinates, chosen uniformly at random, of its vector argument by $d/k$ and sets the other ones to zero. We have  $\omega = \frac{d}{k}-1$. The iteration complexity with $\gamma=\frac{1}{L}$ is $\mathcal{O}\big((\kappa+\frac{d^2}{k^2})\log\frac{1}{\epsilon}\big)$ and the communication complexity, in terms of average number of floats sent by every worker to the master, is $\mathcal{O}\big((k\kappa+\frac{d^2}{k})\log\frac{1}{\epsilon}\big)$, since $k$ floats are sent by every worker at every iteration.  Thus, by choosing $k=\lceil d/\sqrt{\kappa}\rceil$, as long as $d\geq \sqrt{\kappa}$, the communication complexity in terms of floats is $\mathcal{O}\left(d\sqrt{\kappa} \log\frac{1}{\epsilon}\right)$; this is the same as the one of Scaffnew with $\gamma=\frac{1}{L}$ and $p=\frac{1}{\sqrt{\kappa}}$, but  \algn{RandProx-FL} with \texttt{rand}-$k$ compressors removes the necessity to communicate full $d$-dimensional vectors periodically.

\section{Contraction of gradient descent}

\noindent \textbf{Lemma 1.}\ \ \emph{For every $\gamma>0$, the gradient descent operator $\mathrm{Id}-\gamma \nabla \ff$ is $c_{\gamma}$-Lipschitz continuous, with $c_{\gamma}\coloneqq\max(1-\gamma\mu_{\ff},\gamma L_{\ff}-1)$. That is, for every $(x,x')\in\mathcal{X}^2$,
\begin{equation*}
\|(\mathrm{Id}-\gamma \nabla \ff)x-(\mathrm{Id}-\gamma \nabla \ff)x'\|\leq c_{\gamma} \|x-x'\|.
\end{equation*}}

\noindent \textit{Proof}\  \ Let $(x,x')\in\mathcal{X}^2$. By cocoercivity of $\nabla \ff - \mu_{\ff} \mathrm{Id}$, we have \citep[Lemma 3.11]{bub15} $\langle \nabla \ff(x)-\nabla \ff(x') , x-x'\rangle \geq \frac{L_{\ff}\mu_{\ff} }{L_{\ff}+\mu_{\ff}} \|x-x'\|^2 + \frac{1}{L_{\ff}+\mu_{\ff} }  \|\nabla \ff(x)-\nabla \ff(x') \|^2$. Hence,
\begin{align*}
\|(\mathrm{Id}-\gamma \nabla \ff)x-(\mathrm{Id}-\gamma \nabla \ff)x'\|^2&\leq \big( {\textstyle 1-\frac{2\gamma L_{\ff}\mu_{\ff} }{L_{\ff}+\mu_{\ff}}}\big)\|x-x'\|^2 \notag\\
&\quad+  \big({\textstyle \gamma^2- \frac{2\gamma}{L_{\ff}+\mu_{\ff} }}\big)\|\nabla \ff(x)-\nabla \ff(x') \|^2.
\end{align*}
Thus, if $\gamma\leq \frac{2}{L_{\ff}+\mu_{\ff} }$, since $\|\nabla \ff(x)-\nabla \ff(x') \|\geq \mu_{\ff} \|x-x'\|$,
\begin{align*}
\|(\mathrm{Id}-\gamma \nabla \ff)x-(\mathrm{Id}-\gamma \nabla \ff)x'\|^2&\leq \Big( {\textstyle 1-\frac{2\gamma L_{\ff}\mu_{\ff} }{L_{\ff}+\mu_{\ff}}+(\gamma^2- \frac{2\gamma}{L_{\ff}+\mu_{\ff} })\mu_{\ff}^2}\Big)\|x-x'\|^2 \\
&=(1-\gamma\mu_{\ff})^2\|x-x'\|^2.
\end{align*}
On the other hand, if $\gamma\geq \frac{2}{L_{\ff}+\mu_{\ff} }$, since  $\|\nabla \ff(x)-\nabla \ff(x') \|\leq L_{\ff} \|x-x'\|$,
\begin{align*}
\|(\mathrm{Id}-\gamma \nabla \ff)x-(\mathrm{Id}-\gamma \nabla \ff)x'\|^2&\leq \Big( {\textstyle 1-\frac{2\gamma L_{\ff}\mu_{\ff} }{L_{\ff}+\mu_{\ff}}+(\gamma^2- \frac{2\gamma}{L_{\ff}+\mu_{\ff} })L_{\ff}^2}\Big)\|x-x'\|^2 \\
&=(\gamma L_{\ff}-1)^2\|x-x'\|^2.
\end{align*}
Since $\max(1-\gamma\mu_{\ff},\gamma L_{\ff}-1)=( 1-\gamma\mu_{\ff}$ if $\gamma\leq \frac{2}{L_{\ff}+\mu_{\ff} }$, 
$\gamma L_{\ff}-1$ otherwise$)\geq 0$, we arrive at the given expression of $c_{\gamma}$.\hfill$\square$\smallskip

We note that if $\gamma < \frac{2}{L_{\ff}}$ and $\mu_{\ff}>0$, $c_{\gamma}<1$.

\section{Proof of Theorem 1}

Let $t\in\mathbb{N}$. Let $p^t\in\partial \g(\hat{x}^{t})$ be such that $\hat{x}^{t}=x^t-\gamma \nabla \ff(x^t)-\gamma p^t -\gamma K^* u^t$; $p^t$  exists and is unique, by properties of the proximity operator. We also define $p^\star \coloneqq  -\nabla \ff(x^\star) - K^* u^\star$; we have $p^\star\in\partial \g(x^\star)$. Let $q^t \coloneqq  p^t - \mu_{\g} \hat{x}^{t}$ and $q^\star \coloneqq  p^\star - \mu_{\g} x^\star$. We have
$(1+\gamma\mu_{\g})\hat{x}^{t}=x^t-\gamma \nabla \ff(x^t)-\gamma q^t -\gamma K^* u^t$. Let $w^t\coloneqq x^t-\gamma \nabla \ff(x^t)$ and $w^\star\coloneqq x^\star-\gamma \nabla \ff(x^\star)$.

Using  $\hat{u}^{t+1}$ defined in \eqref{eqkolgkr}, we have 
\begin{align*}
\Exp{\sqnorm{x^{t+1}-x^\star}\;|\;\mathcal{F}_t}&=\sqnorm{\Exp{x^{t+1}\;|\;\mathcal{F}_t}-x^\star}+\Exp{\sqnorm{x^{t+1}-\Exp{x^{t+1}\;|\;\mathcal{F}_t}}\;|\;\mathcal{F}_t}\\
&\leq \sqnorm{\hat{x}^{t}-x^\star-\gamma K^*(\hat{u}^{t+1}-u^t)}+\gamma^2\oma\sqnorm{\hat{u}^{t+1}-u^t}\\
&\quad - \gamma^2\zeta\sqnorm{K^*(\hat{u}^{t+1}-u^t)}\ .
\end{align*}
Moreover,
\begin{align*}
\sqnorm{\hat{x}^{t}-x^\star-\gamma K^*(\hat{u}^{t+1}-u^t)}
&=\sqnorm{\hat{x}^{t}-x^\star} +\gamma^2\sqnorm{K^*(\hat{u}^{t+1}-u^t)}\\
&\quad-2\gamma \langle
\hat{x}^{t}-x^\star,K^*(\hat{u}^{t+1}-u^t)\rangle \\
&\leq (1+\gamma\mu_{\g}) \sqnorm{\hat{x}^{t}-x^\star} +\gamma^2\sqnorm{K^*(\hat{u}^{t+1}-u^t)}\\
&\quad-2\gamma \langle
\hat{x}^{t}-x^\star,K^*(\hat{u}^{t+1}-u^\star)\rangle +2\gamma \langle
\hat{x}^{t}-x^\star,K^*(u^t-u^\star)\rangle \\
&=  \langle w^t-w^\star-\gamma (q^t-q^\star) -\gamma K^* (u^t-u^\star),\hat{x}^{t}-x^\star\rangle \\
&\quad+\gamma^2\sqnorm{K^*(\hat{u}^{t+1}-u^t)}\\
&\quad-2\gamma \langle
\hat{x}^{t}-x^\star,K^*(\hat{u}^{t+1}-u^\star)\rangle +2\gamma \langle
\hat{x}^{t}-x^\star,K^*(u^{t}-u^\star)\rangle \\
&=  -2\gamma \langle q^t-q^\star,\hat{x}^{t}-x^\star\rangle\\
&\quad+\langle w^t-w^\star+\gamma (q^t-q^\star) +\gamma K^* (u^t-u^\star),\hat{x}^{t}-x^\star\rangle \\
&\quad+\gamma^2\sqnorm{K^*(\hat{u}^{t+1}-u^t)}-2\gamma \langle
\hat{x}^{t}-x^\star,K^*(\hat{u}^{t+1}-u^\star)\rangle \\
&=  -2\gamma \langle q^t-q^\star,\hat{x}^{t}-x^\star\rangle\\
&\quad+\frac{1}{1+\gamma\mu_{\g}}\langle w^t-w^\star+\gamma (q^t-q^\star) +\gamma K^* (u^t-u^\star),\\
&\quad \qquad w^t-w^\star-\gamma (q^t-q^\star) -\gamma K^* (u^t-u^\star)\rangle \\
&\quad+\gamma^2\sqnorm{K^*(\hat{u}^{t+1}-u^t)}-2\gamma \langle
\hat{x}^{t}-x^\star,K^*(\hat{u}^{t+1}-u^\star)\rangle \\
&=  -2\gamma \langle q^t-q^\star,\hat{x}^{t}-x^\star\rangle+\frac{1}{1+\gamma\mu_{\g}}\sqnorm{w^t-w^\star}\\
&\quad-\frac{\gamma^2}{1+\gamma\mu_{\g}}\sqnorm{q^t-q^\star+K^* (u^t-u^\star)}\\
&\quad+\gamma^2\sqnorm{K^*(\hat{u}^{t+1}-u^t)}-2\gamma \langle
\hat{x}^{t}-x^\star,K^*(\hat{u}^{t+1}-u^\star)\rangle .
\end{align*}
We have $\langle q^t-q^\star,\hat{x}^{t}-x^\star\rangle\geq 0$. Hence, 
\begin{align*}
\sqnorm{\hat{x}^{t}-x^\star-\gamma K^*(\hat{u}^{t+1}-u^t)}
&\leq  \frac{1}{1+\gamma\mu_{\g}}\sqnorm{w^t-w^\star}-\frac{\gamma^2}{1+\gamma\mu_{\g}}\sqnorm{q^t-q^\star+K^* (u^t-u^\star)}\\
&\quad+\gamma^2\sqnorm{K^*(\hat{u}^{t+1}-u^t)}
-2\gamma \langle
\hat{x}^{t}-x^\star,K^*(\hat{u}^{t+1}-u^\star)\rangle,
\end{align*}
so that 
\begin{align*}
\Exp{\sqnorm{x^{t+1}-x^\star}\;|\;\mathcal{F}_t}&\leq  \frac{1}{1+\gamma\mu_{\g}}\sqnorm{w^t-w^\star}-\frac{\gamma^2}{1+\gamma\mu_{\g}}\sqnorm{q^t-q^\star+K^* (u^t-u^\star)}\\
&\quad+\gamma^2(1-\zeta)\sqnorm{K^*(\hat{u}^{t+1}-u^t)}-2\gamma \langle
\hat{x}^{t}-x^\star,K^*(\hat{u}^{t+1}-u^\star)\rangle\\
&\quad+\gamma^2\oma\sqnorm{\hat{u}^{t+1}-u^t}.
\end{align*}

On the other hand, 
\begin{align}
\Exp{\sqnorm{u^{t+1}-u^\star}\;|\;\mathcal{F}_t}&\leq \sqnorm{u^{t}-u^\star+\frac{1}{1+\omega}\big(\hat{u}^{t+1} -u^t\big)}
+\frac{\omega}{(1+\omega)^2}\sqnorm{\hat{u}^{t+1} -u^t }\notag\\
&=\frac{\omega^2}{(1+\omega)^2}\sqnorm{u^{t}-u^\star}+\frac{1}{(1+\omega)^2}\sqnorm{\hat{u}^{t+1}-u^\star}\notag\\
&\quad+\frac{2\omega}{(1+\omega)^2}\langle u^{t}-u^\star,
\hat{u}^{t+1}-u^\star\rangle+\frac{\omega}{(1+\omega)^2}\sqnorm{\hat{u}^{t+1} -u^\star }\notag\\
&\quad+\frac{\omega}{(1+\omega)^2}\sqnorm{u^{t} -u^\star }-\frac{2\omega}{(1+\omega)^2}\langle u^{t}-u^\star,
\hat{u}^{t+1}-u^\star\rangle\notag\\
&=\frac{1}{1+\omega}\sqnorm{\hat{u}^{t+1} -u^\star }+\frac{\omega}{1+\omega}\sqnorm{u^{t} -u^\star }.\label{eq55a}
\end{align}
Let $s^{t+1}\in \partial \hc(\hat{u}^{t+1})$ be such that $\hat{u}^{t+1}=u^t + \tau K\hat{x}^{t}-\tau s^{t+1}$;  $s^{t+1}$ exists and is unique. We also define $s^\star\coloneqq  Kx^\star$; we have $s^\star\in \partial \hc(u^\star)$. 
Therefore,
\begin{align*}
\sqnorm{\hat{u}^{t+1}-u^\star}&=\sqnorm{(u^t-u^\star)+(\hat{u}^{t+1}-u^t)}\\
&=\sqnorm{u^t-u^\star}+\sqnorm{\hat{u}^{t+1}-u^t}+2\langle u^t-u^\star,\hat{u}^{t+1}-u^t\rangle\\
&=\sqnorm{u^t-u^\star}+2 \langle \hat{u}^{t+1}-u^\star,\hat{u}^{t+1}-u^t\rangle - \sqnorm{\hat{u}^{t+1}-u^t}\\
&=\sqnorm{u^t-u^\star} - \sqnorm{\hat{u}^{t+1}-u^t}+2\tau \langle \hat{u}^{t+1}-u^\star, K(\hat{x}^{t}-x^\star)\rangle\\
&\quad -2\tau \langle \hat{u}^{t+1}-u^\star,s^{t+1}-s^\star\rangle.
\end{align*}
Hence, 
\begin{align*}
\frac{1}{\gamma}\Exp{\sqnorm{x^{t+1}-x^\star}\;|\;\mathcal{F}_t}&+\frac{1+\omega}{\tau}\Exp{\sqnorm{u^{t+1}-u^\star}\;|\;\mathcal{F}_t}\\
&\leq  \frac{1}{\gamma(1+\gamma\mu_{\g})}\sqnorm{w^t-w^\star}-\frac{\gamma}{1+\gamma\mu_{\g}}\sqnorm{q^t-q^\star+K^* (u^t-u^\star)}\\
&\quad+\gamma(1-\zeta)\sqnorm{K^*(\hat{u}^{t+1}-u^t)}-2 \langle
\hat{x}^{t}-x^\star,K^*(\hat{u}^{t+1}-u^\star)\rangle\\
&\quad+\gamma\oma\sqnorm{\hat{u}^{t+1}-u^t}+ \frac{1}{\tau}\sqnorm{u^t-u^\star} - \frac{1}{\tau}\sqnorm{\hat{u}^{t+1}-u^t}\\
&\quad+2 \langle \hat{u}^{t+1}-u^\star, K(\hat{x}^{t}-x^\star)\rangle -2 \langle \hat{u}^{t+1}-u^\star,s^{t+1}-s^\star\rangle \\
&\quad+\frac{\omega}{\tau}\sqnorm{u^{t} -u^\star }\\
&\leq \frac{1}{\gamma(1+\gamma\mu_{\g})}\sqnorm{w^t-w^\star}-\frac{\gamma}{1+\gamma\mu_{\g}}\sqnorm{q^t-q^\star+K^* (u^t-u^\star)}\\
&\quad+\frac{1+\omega}{\tau}\sqnorm{u^{t} -u^\star }+ \left(\gamma \big((1-\zeta)\|K\|^2+\oma\big) - \frac{1}{\tau}\right) \sqnorm{\hat{u}^{t+1}-u^t} \\
&\quad-2 \langle \hat{u}^{t+1}-u^\star,s^{t+1}-s^\star\rangle\\
&\leq \frac{1}{\gamma(1+\gamma\mu_{\g})}\sqnorm{w^t-w^\star}-\frac{\gamma}{1+\gamma\mu_{\g}}\sqnorm{q^t-q^\star+K^* (u^t-u^\star)}\\
&\quad+\frac{1+\omega}{\tau}\sqnorm{u^{t} -u^\star }-2 \langle \hat{u}^{t+1}-u^\star,s^{t+1}-s^\star\rangle.
\end{align*}
By $\mu_{\hc}$-strong monotonicity of $\partial \hc$, $\langle \hat{u}^{t+1}-u^\star,s^{t+1}-s^\star\rangle \geq \mu_{\hc} \sqnorm{\hat{u}^{t+1}-u^\star}$, and using \eqref{eq55a},
 \begin{equation*}
 \langle \hat{u}^{t+1}-u^\star,s^{t+1}-s^\star\rangle \geq \mu_{\hc} \left((1+\omega)\Exp{\sqnorm{u^{t+1}-u^\star}\;|\;\mathcal{F}_t}-\omega \sqnorm{u^{t} -u^\star }\right).
  \end{equation*}
Hence, 
\begin{align}
\frac{1}{\gamma}\Exp{\sqnorm{x^{t+1}-x^\star}\;|\;\mathcal{F}_t}&+(1+\omega)\left(\frac{1}{\tau}+2\mu_{\hc}\right)\Exp{\sqnorm{u^{t+1}-u^\star}\;|\;\mathcal{F}_t}\notag\\
& \leq \frac{1}{\gamma(1+\gamma\mu_{\g})}\sqnorm{w^t-w^\star}-\frac{\gamma}{1+\gamma\mu_{\g}}\sqnorm{q^t-q^\star+K^* (u^t-u^\star)}\notag\\
&\quad+\left(\frac{1+\omega}{\tau}+2\omega \mu_{\hc} \right)
\sqnorm{u^{t} -u^\star }.\label{gkehgka}
\end{align}
After Lemma 1, 
\begin{eqnarray*}
\sqnorm{w^t-w^\star}&=& \sqnorm{(\mathrm{Id}-\gamma\nabla \ff)x^t-(\mathrm{Id}-\gamma\nabla \ff)x^\star} \\
&\leq &\max(1-\gamma\mu_{\ff},\gamma L_{\ff}-1)^2 \sqnorm{x^t-x^\star}.
\end{eqnarray*}
Plugging this inequality  in \eqref{gkehgka} yields 
\begin{align}
\Exp{\Psi^{t+1}\;|\;\mathcal{F}_t}
&\leq \frac{1}{\gamma(1+\gamma\mu_{\g})}\max(1-\gamma\mu_{\ff},\gamma L_{\ff}-1)^2 \sqnorm{x^t-x^\star}\label{eqrec2a}\\
&\quad+\left(\frac{1+\omega}{\tau}+2\omega \mu_{\hc}\right)\sqnorm{u^{t} -u^\star }-\frac{\gamma}{1+\gamma\mu_{\g}}\sqnorm{q^t-q^\star+K^* (u^t-u^\star)}.\notag
\end{align}
Ignoring the last term in \eqref{eqrec2a}, we obtain: 
\begin{align}
\Exp{\Psi^{t+1}\;|\;\mathcal{F}_t}
&\leq \max\left(\frac{(1-\gamma\mu_{\ff})^2}{1+\gamma\mu_{\g}},\frac{(\gamma L_{\ff}-1)^2}{1+\gamma\mu_{\g}},1-\frac{2\tau \mu_{\hc} }{(1+\omega)(1+2\tau \mu_{\hc})}\right)\Psi^t.\label{eqrec2b}
\end{align}
Using the tower rule, we can unroll the recursion in \eqref{eqrec2b} to obtain the unconditional expectation of $\Psi^{t+1}$. Since $\Exp{\Psi^t}\rightarrow 0$, we have $\Exp{\sqnorm{x^t-x^\star}}\rightarrow 0$ and $\Exp{\sqnorm{u^t-u^\star}}\rightarrow 0$. Moreover, using classical results on supermartingale convergence \citep[Proposition A.4.5]{ber15}, it follows from \eqref{eqrec2b} that $\Psi^t\rightarrow 0$ almost surely. Almost sure convergence of $x^t$ and $u^t$ follows. Finally, 
by Lipschitz continuity of $\nabla \ff$, $K^*$, $\mathrm{prox}_{\g}$, we can upper bound $\|\hat{x}^t-x^\star\|^2$ by a linear combination of $\|x^t-x^\star\|^2$ and $\|u^t-u^\star\|^2$. 
It follows that $\Exp{\sqnorm{\hat{x}^t-x^\star}}\rightarrow 0$ linearly with the same rate $c$ and that $\hat{x}^t \rightarrow x^\star$ almost surely, as well.
\hfill$\square$

\section{Proof of Theorem 2}

Let us go back to \eqref{eqrec2a}. Since $\g=0$, we have $q^t=q^\star=0$ and $\mu_{\g}=0$, so that 
\begin{align*}
\Exp{\Psi^{t+1}\;|\;\mathcal{F}_t}
&\leq \frac{1}{\gamma}\max(1-\gamma\mu_{\ff},\gamma L_{\ff}-1)^2 \sqnorm{x^t-x^\star}\notag+\left(\frac{1+\omega}{\tau}+2\omega \mu_{\hc}\right)\sqnorm{u^{t} -u^\star }\notag\\
&\quad -\gamma\sqnorm{K^* (u^t-u^\star)}.
\end{align*}
We have $\sqnorm{K^*(u^t-u^\star)}\geq \lambda_{\min}(KK^*)\sqnorm{u^t-u^\star}$.
This yields
\begin{align}
\Exp{\Psi^{t+1}\;|\;\mathcal{F}_t}
&\leq \frac{1}{\gamma}\max(1-\gamma\mu_{\ff},\gamma L_{\ff}-1)^2 \sqnorm{x^t-x^\star}\notag\\
&\quad+\left(\frac{1+\omega}{\tau}+2\omega \mu_{\hc}- \gamma\lambda_{\min}(KK^*) \right)\sqnorm{u^{t} -u^\star }\notag\\
&\leq \max\left((1-\gamma\mu_{\ff})^2,(\gamma L_{\ff}-1)^2,1-\frac{2\tau \mu_{\hc} + \gamma\tau\lambda_{\min}(KK^*)}{(1+\omega)(1+2\tau \mu_{\hc})}\right)\Psi^t.\label{eqrec2d}
\end{align}
The end of the proof is the same as the one of Theorem 1.
\hfill$\square$\medskip

Let us add here a remark on the PAPC algorithm, which is the particular case of \algn{RandProx} when $\omega=0$, in the conditions of Theorem 2:\medskip

\noindent \textbf{Remark 2} (PAPC vs.\ proximal gradient descent on the dual problem)\ \ If  $\mu_{\ff}>0$, $\ffc$ is $\mu^{-1}$-smooth and $L_{\ff}^{-1}$-strongly convex. Then $\ffc \circ -K^*$ is $\mu_{\ff}^{-1}\|K\|^2$-smooth and $L_{\ff}^{-1}\lambda_{\min}(KK^*)$-strongly convex. 
So, if $\nabla \ffc$ is computable,  one can apply the proximal gradient algorithm on the dual problem \eqref{eqpbd}, which iterates 
$u^{t+1}=\mathrm{prox}_{\tau \hc}\big(u^{t}+\tau K\nabla \ffc(-K^*u^{t})\big)$, with $\tau\in \big(0,\frac{2\mu_{\ff}}{ \|K\|^2}\big)$.  If
$\lambda_{\min}(KK^*)>0$,
this algorithm converges linearly: $\|u^{t+1}-u^\star\|^2\leq c^2 \|u^{t}-u^\star\|^2$ with $c = \max\big(1-\tau L_{\ff}^{-1}\lambda_{\min}(KK^*),\tau\mu_{\ff}^{-1}\|K\|^2-1\big)$. $c$ is smallest with $\tau=2/\big(\mu_{\ff}^{-1} \|K\|^2+L_{\ff}^{-1}\lambda_{\min}(KK^*)\big)$, in which case
\begin{equation*}
c=\frac{1-\frac{\mu_{\ff}}{L_{\ff}}\frac{\lambda_{\min}(KK^*)}{\|K\|^2}}{1+\frac{\mu_{\ff}}{L_{\ff}}\frac{\lambda_{\min}(KK^*)}{\|K\|^2}}.
\end{equation*}
This is much worse than the rate of the PAPC algorithm, since it involves the product of the condition numbers $L_{\ff}/\mu_{\ff}$ and $\|K\|^2/\lambda_{\min}(KK^*)$, instead of their maximum. This is due to calling gradients of $\ffc \circ -K^*$, whereas $\ff$ and $K$ are split, or decoupled, in the PAPC algorithm.

\section{Proof of Theorem 4 and further discussion}

\begin{figure*}[t]
\centering
	\begin{minipage}{.5\textwidth}
	\begin{algorithm}[H]
		\caption{\algn{RandPriLiCo} [new]}
		\begin{algorithmic}
			\STATE  \textbf{input:} initial points $x^0\in\mathcal{X}$, $v^0\in \mathrm{ran}(W)$; 
			\STATE stepsizes $\gamma>0$, $\tau>0$; $\omega\geq 0$
			\FOR{$t=0, 1, \ldots$}
			\STATE $\hat{x}^{t} \coloneqq  x^t -\gamma \nabla \ff(x^t) - \gamma v^t$
			\STATE $d^{t+1}\coloneqq  \tau\mathcal{S}^t(  W \hat{x}^{t}- a)$
			\STATE $v^{t+1}\coloneqq v^{t}+\frac{1}{1+\omega}d^{t+1}$
			\STATE  $x^{t+1} \coloneqq \hat{x}^{t}-\gamma d^{t+1}$
			\ENDFOR
		\end{algorithmic}
	\end{algorithm}\end{minipage}\end{figure*}

We observe that in \algn{RandProx-LC} and 
Theorem 4, it is as if the sequence $(u^t_0)_{t\in\mathbb{N}}$ had been computed by the following iteration, initialized with $x^0\in\mathcal{X}$ and $u^0_0\coloneqq P_{\mathrm{ran}(K)}(u^0)$:
 \begin{equation*}
 \left\lfloor
 \begin{array}{l}
\hat{x}^{t} \coloneqq  x^t -\gamma \nabla \ff(x^t) - \gamma v^t\\
u_0^{t+1}\coloneqq u_0^t +\frac{1}{1+\omega}P_{\mathrm{ran}(K)}\Rop^t\big(\tau (K \hat{x}^{t}-b)\big)\\
v^{t+1}\coloneqq K^* u_0^{t+1}\\
x^{t+1} \coloneqq \hat{x}^{t}-\gamma (1+\omega) (v^{t+1}-v^t)
\end{array}\right..
 \end{equation*}
Then we remark that this is simply the iteration of \algn{RandProx}, with $\Rop^t$ replaced by $\Ropp^t\coloneqq P_{\mathrm{ran}(K)}\Rop^t$. Since its argument $r^t = \tau (K \hat{x}^{t}-b)$ is always in $\mathrm{ran}(K)$, $\Ropp^t$ is unbiased, and we have, for every $t\geq 0$, 
 \begin{equation*}
\Exp{\sqnorm{\Ropp^t(r^t)-r^t}\;|\;\widetilde{\mathcal{F}}_t}\leq \Exp{\sqnorm{\Rop^t(r^t)-r^t}\;|\;\widetilde{\mathcal{F}}_t}\leq \omega \sqnorm{r^t},
 \end{equation*}
 where $\widetilde{\mathcal{F}}_t$ the $\sigma$-algebra generated by the collection of random variables $(x^0,u_0^0),\ldots, (x^t,u_0^t)$. 
 Also, $\oma$ is unchanged.  
Therefore, the analysis of \algn{RandProx} in Theorem 2 applies, with $u^t$ replaced by $u_0^t$ and $u^\star$ by $u_0^\star$. Now, for every $u\in \mathrm{ran}(K)$, 
 \begin{equation*}
\sqnorm{K^*u}\geq \lambda_{\min}^+(KK^*)\sqnorm{u},
\end{equation*}
and using this lower bound in the proof of Theorem 2, with $\mu_{\hc}=0$, we obtain Theorem 4.\hfill$\square$\medskip

Furthermore, the constraint $Kx=b$ is equivalent to the constraint $K^*Kx=K^*b$; so, let us consider problems where we are given $K^*K$ and not $K$ in the first place:

Let $W$ be a  linear operator on $\mathcal{X}$, which is self-adjoint, i.e.\ $W^*=W$, and positive, i.e.\ $\langle Wx,x\rangle \geq 0$ for every $x\in\mathcal{X}$. 
Let $a\in\mathrm{ran}(W)$. We consider the linearly constrained minimization problem
\begin{equation}
\mathrm{Find} \ x^\star \in \argmin_{x\in\mathcal{X}} \, \ff(x) \quad \mbox{s.t.} \quad Wx=a.\label{eqlck1w}
 \end{equation}
 Now, we let $\mathcal{U}\coloneqq \mathcal{X}$ and $K=K^*\coloneqq\sqrt{W}$, where $\sqrt{W}$ is the unique positive self-adjoint linear operator on $\mathcal{X}$ such that $\sqrt{W}\sqrt{W}=W$. Also, $b$ is defined as the unique element in $\mathrm{ran}(W)=\mathrm{ran}(K)$ such that $\sqrt{W}b=a$. Then \eqref{eqlck1w} is equivalent to \eqref{eqlck1} and the dual problem is \eqref{eqpbd2}.  We consider the Randomized Primal Linearly Constrained minimization algorithm (\algn{RandPriLiCo}), shown above.  We suppose that the stochastic operators $\mathcal{S}^t$ in \algn{RandPriLiCo} satisfy, for every $t\geq 0$,
 \begin{equation}
\Exp{\mathcal{S}^t(r^t)\;|\;\widetilde{\mathcal{F}}_t}=r^t\quad\mbox{and}\quad\Exp{\sqnorm{\mathcal{S}^t(r^t)-r^t}\;|\;\widetilde{\mathcal{F}}_t} \leq \omega \sqnorm{r^t},
\end{equation}
 for some $\omega\geq 0$, where $r^t \coloneqq \tau W \hat{x}^{t}-\tau a$.
 
 In addition, we suppose that the $\mathcal{S}^t$ commute with $\sqrt{W}$: for every $t\geq 0$ and $x\in \mathcal{X}$,
 \begin{equation*}
\sqrt{W}\mathcal{S}^t(x) = \mathcal{S}^t(\sqrt{W}x).
 \end{equation*}
 This is satisfied with the Bernoulli operators or some linear sketching operators, for instance. Then \algn{RandPriLiCo} is equivalent to \algn{RandProx-LC}, with $\mathcal{S}^t$ playing the role of $\Rop^t$ and $\oma=\|W\|\omega$, $\zeta=0$. Applying Theorem 4 with these equivalences, we obtain:\medskip
 
 \noindent\textbf{Theorem 6.}\ \  \emph{In the setting of \eqref{eqlck1w}, suppose that $\mu_{\ff}>0$. In \algn{RandPriLiCo}, suppose that $0<\gamma < \frac{2}{L_{\ff}}$, $\tau>0$ and $\gamma \tau \|W\|(1+\omega)\leq 1$. 
Define the Lyapunov function, for every $t\geq 0$,
\begin{equation}
\Psi^{t}\coloneqq  \frac{1}{\gamma}\sqnorm{x^{t}-x^\star}+\frac{1+\omega}{\tau}\sqnorm{u^{t}_0-u_0^\star},
\end{equation}
where $u^t_0$ is the unique element in $\mathrm{ran}(W)$ such that $v^t=\sqrt{W} u^t_0$, $x^\star$ is the unique solution of \eqref{eqlck1w} and $u_0^\star$ is the unique element in $\mathrm{ran}(W)$ such that $-\nabla \ff(x^\star)=\sqrt{W}u_0^\star$. Then \algn{RandPriLiCo} converges linearly:  for every $t\geq 0$,
\begin{align}
\Exp{\Psi^{t}}&\leq c^t \Psi^0,
\end{align}
where 
\begin{equation}
c\coloneqq \max\left((1-\gamma\mu_{\ff})^2,(\gamma L_{\ff}-1)^2,1-\frac{\gamma\tau\lambda_{\min}^+(W)}{1+\omega}\right)<1.
\end{equation}
Also, $(x^t)_{t\in\mathbb{N}}$ and $(\hat{x}^t)_{t\in\mathbb{N}}$ both converge  to $x^\star$ almost surely.
}\medskip

 \algn{RandPriLiCo}  can be applied to decentralized optimization, like in \citet{kov20,sal22} but with randomized communication; we leave the detailed study of this setting for future work.

\section{Particular case $\ff=0$: randomized Chambolle--Pock algorithm}\label{seccp}

\begin{figure*}[t]
\begin{minipage}{.48\textwidth}
\begin{algorithm}[H]
		\caption{CP algorithm\\ \citep{cha11a}}
		\begin{algorithmic}
			\STATE  \textbf{input:} initial points $x^0\in\mathcal{X}$, $u^0\in\mathcal{U}$; 
			\STATE stepsizes $\gamma>0$, $\tau>0$
			\STATE $\hat{x}^{0} \coloneqq  \mathrm{prox}_{\gamma \g}\big(x^{0} - \gamma K^* u^{0}\big)$
			\FOR{$t=0, 1, \ldots$}	
			\STATE $u^{t+1}\coloneqq \mathrm{prox}_{\tau \hc} \big(u^t+\tau K \hat{x}^{t}\big)$
			\STATE  // $x^{t+1} \coloneqq \hat{x}^{t}-\gamma K^*(u^{t+1}-u^t)$
			\STATE $\hat{x}^{t+1} \coloneqq  \mathrm{prox}_{\gamma \g}\big(\hat{x}^{t}- \gamma K^*(2u^{t+1}-u^t)\big)$
			\ENDFOR
		\end{algorithmic}
	\end{algorithm}\end{minipage}\ \ \ \ \ \begin{minipage}{.48\textwidth}
	\begin{algorithm}[H]
		\caption{\algn{RandProx-CP} [new]}
		\begin{algorithmic}
			\STATE  \textbf{input:} initial points $x^0\in\mathcal{X}$, $u^0\in\mathcal{U}$; 
			\STATE stepsizes $\gamma>0$, $\tau>0$; $\omega\geq 0$
			\STATE $\hat{x}^{0} \coloneqq  \mathrm{prox}_{\gamma \g}\big(x^{0} - \gamma K^* u^{0}\big)$
			\FOR{$t=0, 1, \ldots$}
			\STATE $d^t\coloneqq \Rop^t\big(\mathrm{prox}_{\tau \hc} (u^t+\tau K \hat{x}^{t})-u^t\big)$
			\STATE $u^{t+1}\coloneqq u^t +\frac{1}{1+\omega}d^t$
			 \STATE //  $x^{t+1} \coloneqq \hat{x}^{t}-\gamma K^*d^t$
			\STATE $\hat{x}^{t+1} \coloneqq  \mathrm{prox}_{\gamma \g}\big(\hat{x}^{t}- \gamma K^*(u^{t+1}+d^t)\big)$
			\ENDFOR
		\end{algorithmic}
	\end{algorithm}\end{minipage}\end{figure*}

In this section, we suppose that $\ff=0$. 
The primal problem \eqref{eqpb0} becomes:
\begin{equation}
\mathrm{Find} \ x^\star \in \argmin_{x\in\mathcal{X}}  \Big( \g(x) + \h(Kx)\Big),\label{eqpb0ab}
 \end{equation}
 and the dual problem \eqref{eqpbd} becomes:
\begin{equation}
\mathrm{Find} \ u^\star \in \argmin_{u\in\mathcal{U}} \, \Big( \gc(-K^*u)+\hc(u)\Big).\label{eqpbdab}
 \end{equation}
 The PDDY algorithm becomes the Chambolle-Pock (CP), a.k.a.\ PDHG, algorithm \citep{cha11a}, shown above. \algn{RandProx} can be rewritten as \algn{RandProx-CP}, shown above, too. In both algorithms, the variable $x^t$ is not needed any more and can be removed. 
 
 Since $\ff=0$, $L_{\ff}>0$ can be set arbitrarily close to zero, so that Theorem 1 can be rewritten as:\medskip
 
\noindent\textbf{Theorem 7.}\ \  \emph{Suppose that $\mu_{\g}>0$ and $\mu_{\hc}>0$.
In \algn{RandProx-CP}, suppose that $\gamma>0$, $\tau>0$, $\gamma \tau \big((1-\zeta)\|K\|^2+\oma\big)\leq 1$. 
Define the Lyapunov function, for every $t\geq 0$,
\begin{equation}
\Psi^{t}\coloneqq  \frac{1}{\gamma}\sqnorm{x^{t}-x^\star}+(1+\omega)\left(\frac{1}{\tau}+2\mu_{\hc}\right)\sqnorm{u^{t}-u^\star},
\end{equation}
where $x^\star$ and $u^\star$ are the unique solutions to \eqref{eqpb0ab} and \eqref{eqpbdab}, respectively.
Then \algn{RandProx-CP} converges linearly:  for every $t\geq 0$,
\begin{align}
\Exp{\Psi^{t}}&\leq c^t \Psi^0,
\end{align}
where 
\begin{align}
c&\coloneqq  \max\left(\frac{1}{1+\gamma\mu_{\g}},1-\frac{2\tau \mu_{\hc}}{(1+\omega)(1+2\tau \mu_{\hc})}\right)\\
&\;=1-\min\left(\frac{\gamma\mu_{\g}}{1+\gamma\mu_{\g}},\frac{2\tau \mu_{\hc} }{(1+\omega)(1+2\tau \mu_{\hc})}\right)<1.
\end{align}
Also, $(x^t)_{t\in\mathbb{N}}$ and $(\hat{x}^t)_{t\in\mathbb{N}}$ both converge  to $x^\star$ and $(u^t)_{t\in\mathbb{N}}$ converges to $u^\star$, almost surely.
}\medskip

It would be interesting to study whether the mechanism in the stochastic PDHG algorithm proposed in \citet{cha18} can be viewed as a particular case of \algn{RandProx-CP}; we leave the analysis of this connection for future work. In any case, the strong convexity constants $\mu_{\g}$ and $\mu_{\hc}$ need to be known in the linearly converging version of the stochastic PDHG algorithm, which is not the case here; this is an important advantage of \algn{RandProx-CP}.\smallskip

\begin{figure*}[t]
\begin{minipage}{.42\textwidth}
\begin{algorithm}[H]
		\caption{ADMM}
		\begin{algorithmic}
			\STATE  \textbf{input:} initial points $x^0\in\mathcal{X}$, $u^0\in\mathcal{U}$;
			\STATE stepsize $\gamma>0$
			\FOR{$t=0, 1, \ldots$}	
			\STATE $\hat{x}^{t} \coloneqq  \mathrm{prox}_{\gamma \g}(x^{t}-\gamma u^{t})$
			\STATE $x^{t+1}\coloneqq \mathrm{prox}_{\gamma \h} (\hat{x}^{t}+\gamma u^t)$
			\STATE $u^{t+1}\coloneqq u^t + \frac{1}{\gamma}(\hat{x}^{t}-x^{t+1})$			
			\ENDFOR
		\end{algorithmic}
	\end{algorithm}\end{minipage}\ \ \ \ \ \begin{minipage}{.54\textwidth}
	\begin{algorithm}[H]
		\caption{\algn{RandProx-ADMM} [new]}
		\begin{algorithmic}
			\STATE  \textbf{input:} initial points $x^0\in\mathcal{X}$, $u^0\in\mathcal{U}$; 
			\STATE stepsize $\gamma>0$; $\omega\geq 0$
			\FOR{$t=0, 1, \ldots$}
			\STATE $\hat{x}^{t} \coloneqq  \mathrm{prox}_{\gamma \g}\big(x^{t}- \gamma u^{t}\big)$
			\STATE $d^t\coloneqq \Rop^t\big(  \hat{x}^{t}-\mathrm{prox}_{\gamma(1+\omega) \h} (\hat{x}^{t}+\gamma(1+\omega)u^t)\big)$
			 \STATE $x^{t+1} \coloneqq \hat{x}^{t}-\frac{1}{1+\omega}d^t$
			\STATE $u^{t+1}\coloneqq u^t +\frac{1}{\gamma(1+\omega)^2}d^t$
			\ENDFOR
		\end{algorithmic}
	\end{algorithm}\end{minipage}\end{figure*}

Now, let us look at the particular case $K=\mathrm{Id}$ in \eqref{eqpb0ab} and \eqref{eqpbdab}. 
The primal problem  becomes:
\begin{equation}
\mathrm{Find} \ x^\star \in \argmin_{x\in\mathcal{X}}  \Big( \g(x) + \h(x)\Big),\label{eqpb0abk}
 \end{equation}
 and the dual problem  becomes:
\begin{equation}
\mathrm{Find} \ u^\star \in \argmin_{u\in\mathcal{U}} \, \Big( \gc(-u)+\hc(u)\Big).\label{eqpbdabk}
 \end{equation}
When $K=\mathrm{Id}$, the CP algorithm with $\tau=\frac{1}{\gamma}$ reverts to the Douglas--Rachford algorithm, which is equivalent to the Alternating Direction Method of Multipliers (ADMM) \citep{boy11,con19}, shown above. Therefore, in that case, with $\oma=\omega$, $\zeta=0$ and $\tau=\frac{1}{\gamma(1+\omega)}$, \algn{RandProx-CP} can be rewritten as  \algn{RandProx-ADMM}, shown above. Theorem 7 becomes:\medskip

\noindent\textbf{Theorem 8.}\ \  \emph{Suppose that $\mu_{\g}>0$ and $\mu_{\hc}>0$.
In \algn{RandProx-ADMM}, suppose that $\gamma>0$. 
For every $t\geq 0$, define the Lyapunov function 
\begin{equation}
\Psi^{t}\coloneqq  \frac{1}{\gamma}\sqnorm{x^{t}-x^\star}+(1+\omega)\big(\gamma(1+\omega)+2\mu_{\hc}\big)\sqnorm{u^{t}-u^\star},
\end{equation}
where $x^\star$ and $u^\star$ are the unique solutions to \eqref{eqpb0abk} and \eqref{eqpbdabk}, respectively.
Then \algn{RandProx-ADMM} converges linearly:  for every $t\geq 0$,
\begin{align}
\Exp{\Psi^{t}}&\leq c^t \Psi^0,
\end{align}
where 
\begin{align}
c&\coloneqq  \max\left(\frac{1}{1+\gamma\mu_{\g}},1-\frac{2\tau \mu_{\hc}}{(1+\omega)(1+2\tau \mu_{\hc})}\right)\\
&\;=1-\min\left(\frac{\gamma\mu_{\g}}{1+\gamma\mu_{\g}},\frac{2\tau \mu_{\hc} }{(1+\omega)(1+2\tau \mu_{\hc})}\right)<1.
\end{align}
Also, $(x^t)_{t\in\mathbb{N}}$ and $(\hat{x}^t)_{t\in\mathbb{N}}$ both converge  to $x^\star$ and $(u^t)_{t\in\mathbb{N}}$ converges to $u^\star$, almost surely.
}

\section{Particular case $K=\mathrm{Id}$: randomized Davis--Yin algorithm}\label{seccdy}

\begin{figure*}[t]
\begin{minipage}{.44\textwidth}
\begin{algorithm}[H]
		\caption{DY algorithm\\ \citep{dav17}}
		\begin{algorithmic}
			\STATE  \textbf{input:} initial points $x^0\in\mathcal{X}$, $u^0\in\mathcal{X}$; 
			\STATE stepsize $\gamma>0$
			\FOR{$t=0, 1, \ldots$}
			\STATE $\hat{x}^{t} \coloneqq  \mathrm{prox}_{\gamma \g}\big(x^t -\gamma \nabla \ff(x^t) - \gamma u^t\big)$%
			\STATE $x^{t+1}\coloneqq \mathrm{prox}_{\gamma \h} (\hat{x}^{t}+\gamma u^t)$ 
			\STATE $u^{t+1}\coloneqq u^t + \frac{1}{\gamma}(\hat{x}^{t}-x^{t+1})$
			\ENDFOR
		\end{algorithmic}
	\end{algorithm}\end{minipage}\ \ \ \ \ \begin{minipage}{.52\textwidth}
	\begin{algorithm}[H]
		\caption{\algn{RandProx-DY} [new]}
		\begin{algorithmic}
			\STATE  \textbf{input:} initial points $x^0\in\mathcal{X}$, $u^0\in\mathcal{X}$;   
			\STATE stepsize $\gamma>0$; $\omega\geq 0$
			\FOR{$t=0, 1, \ldots$}
			\STATE $\hat{x}^{t} \coloneqq  \mathrm{prox}_{\gamma \g}\big(x^t -\gamma \nabla \ff(x^t) - \gamma u^t\big)$
			\STATE $d^t\coloneqq \Rop^t\big(  \hat{x}^{t}-\mathrm{prox}_{\gamma(1+\omega) \h} (\hat{x}^{t}+\gamma(1+\omega)u^t)\big)$
			 \STATE $x^{t+1} \coloneqq \hat{x}^{t}-\frac{1}{1+\omega}d^t$
			\STATE $u^{t+1}\coloneqq u^t +\frac{1}{\gamma(1+\omega)^2}d^t$
			\ENDFOR
		\end{algorithmic}
	\end{algorithm}\end{minipage}\end{figure*}

After the particular case $\g=0$ discussed in Section~\ref{secpapc} and the particular case $\ff=0$ discussed in Section~\ref{seccp}, we discuss in this section the third particular case $K=\mathrm{Id}$ in \eqref{eqpb0} and \eqref{eqpbd}. 
The primal problem  becomes:
\begin{equation}
\mathrm{Find} \ x^\star \in \argmin_{x\in\mathcal{X}}  \Big( \ff(x) + \g(x) + \h(x)\Big),\label{eqpb0dy}
 \end{equation}
 and the dual problem  becomes:
\begin{equation}
\mathrm{Find} \ u^\star \in \argmin_{u\in\mathcal{U}} \, \Big( (\ff+\g)^*(-u)+\hc(u)\Big).\label{eqpbddy}
 \end{equation}
 When $K=\mathrm{Id}$, the PDDY algorithm with $\tau=\frac{1}{\gamma}$ reverts to the Davis--Yin (DY) algorithm \citep{dav17}, shown above. 
 Therefore, in that case, with $\oma=\omega$, $\zeta=0$ and $\tau=\frac{1}{\gamma(1+\omega)}$, \algn{RandProx} can be rewritten as  \algn{RandProx-DY}, shown above, too. When $\g=0$,  \algn{RandProx-DY} reverts to \algn{RandProx-FB} and when  $\ff=0$,  \algn{RandProx-DY} reverts to \algn{RandProx-ADMM}; in other words, \algn{RandProx-DY} generalizes \algn{RandProx-FB} and  \algn{RandProx-ADMM} into a single algorithm. Theorem 1 yields:\medskip
 
 \noindent\textbf{Theorem 9.}\ \  \emph{Suppose that $\mu_{\ff}>0$ or $\mu_{\g}>0$, and that $\mu_{\hc}>0$.
In \algn{RandProx-DY}, suppose that $0<\gamma < \frac{2}{L_{\ff}}$. 
For every $t\geq 0$, define the Lyapunov function,
\begin{equation}
\Psi^{t}\coloneqq  \frac{1}{\gamma}\sqnorm{x^{t}-x^\star}+(1+\omega)\big(\gamma(1+\omega)+2\mu_{\hc}\big)\sqnorm{u^{t}-u^\star},
\end{equation}
where $x^\star$ and $u^\star$ are the unique solutions to \eqref{eqpb0dy} and \eqref{eqpbddy}, respectively.
Then  \algn{RandProx-DY} converges linearly:  for every $t\geq 0$,
\begin{align}
\Exp{\Psi^{t}}&\leq c^t \Psi^0,
\end{align}
where 
\begin{equation}
c\coloneqq  \max\left(\frac{(1-\gamma\mu_{\ff})^2}{1+\gamma\mu_{\g}},\frac{(\gamma L_{\ff}-1)^2}{1+\gamma\mu_{\g}},1-\frac{\frac{2}{\gamma} \mu_{\hc}}{(1+\omega)\big(1+\omega+\frac{2}{\gamma} \mu_{\hc}\big)}\right)<1.
\end{equation}
Also, $(x^t)_{t\in\mathbb{N}}$ and $(\hat{x}^t)_{t\in\mathbb{N}}$ both converge  to $x^\star$ and $(u^t)_{t\in\mathbb{N}}$ converges to $u^\star$, almost surely.
}\medskip

We note that in Theorem 9, $\mu_{\hc}>0$ is required. It is only in the case $\g=0$, when  \algn{RandProx-DY} reverts to \algn{RandProx-FB}, that one can apply Theorem 3, which does not require strong convexity of $\hc$.

\section{Proof of Theorem 11}

\noindent\textit{Proof of Theorem 11}\ \   
We have, for every $(x,x')\in\mathcal{X}^2$,
\begin{align*}
\|(\mathrm{Id}-\gamma \nabla \ff)x-(\mathrm{Id}-\gamma \nabla \ff)x'\|^2&= \|x-x'\|^2-  2\gamma\langle \nabla \ff(x)-\nabla \ff(x'),x-x'\rangle \\
&\quad+ \gamma^2\|\nabla \ff(x)-\nabla \ff(x') \|^2\\
&\leq \|x-x'\|^2-  (2\gamma-\gamma^2 L_\ff)\langle \nabla \ff(x)-\nabla \ff(x'),x-x'\rangle,
\end{align*}
where the second inequality follows from cocoercivity of the gradient. 
Moreover, for every $(x,x')\in\mathcal{X}^2$, $D_\ff(x,x')\leq \langle \nabla \ff(x)-\nabla \ff(x'),x-x'\rangle$. 
Therefore, in the proof of Theorem 1, for every primal--dual solution $(x^\star,u^\star)$ and $t\geq 0$, since $\sqnorm{w^t-w^\star}=\sqnorm{(\mathrm{Id}-\gamma\nabla \ff)x^t-(\mathrm{Id}-\gamma\nabla \ff)x^\star}$, 
 \eqref{gkehgka} yields 
\begin{align*}
\Exp{\Psi^{t+1}\;|\;\mathcal{F}_t}
&\leq \frac{1}{\gamma} \sqnorm{x^t-x^\star}-  (2-\gamma L_\ff)D_\ff(x^t,x^\star)\\
&\quad+\left(\frac{1+\omega}{\tau}+2\omega \mu_{\hc}\right)\sqnorm{u^{t} -u^\star }-\gamma\sqnorm{q^t-q^\star+K^* (u^t-u^\star)}.\notag
\end{align*}
Ignoring the last term, this yields
\begin{align}
\Exp{\Psi^{t+1}\;|\;\mathcal{F}_t}
&\leq \frac{1}{\gamma} \sqnorm{x^t-x^\star}+c(1+\omega)\left(\frac{1}{\tau}+2\mu_{\hc}\right)\sqnorm{u^{t}-u^\star}\label{eqsm2}\\
&\quad- (2-\gamma L_\ff)D_\ff(x^t,x^\star)\notag\\
&\leq \Psi^t-  (2-\gamma L_\ff)D_\ff(x^t,x^\star),\label{eqsm1}
\end{align}
with $c=1-\frac{2\tau \mu_{\hc} }{(1+\omega)(1+2\tau \mu_{\hc})}$ in \eqref{eqsm2}. Using classical results on supermartingale convergence \citep[Proposition A.4.5]{ber15}, it follows from \eqref{eqsm1} that $\Psi^t$ converges almost surely to a random variable $\Psi^\infty$ and that
\begin{equation*}
\sum_{t=0}^\infty D_\ff(x^t,x^\star) < +\infty\quad\mbox{almost surely}.
\end{equation*}
Hence, $D_\ff(x^t,x^\star)\rightarrow 0$ almost surely. Moreover, for every $T\geq 0$, 
\begin{equation}
(2-\gamma L_\ff) \sum_{t=0}^T \Exp{D_\ff(x^t,x^\star)} \leq \Psi^0-\Exp{\Psi^{T+1}}\leq \Psi^0\label{eqjghhrg}
\end{equation}
and
\begin{equation*}
(2-\gamma L_\ff) \sum_{t=0}^\infty \Exp{D_\ff(x^t,x^\star)}\leq \Psi^0.
\end{equation*}
Therefore, $\Exp{D_\ff(x^t,x^\star)}\rightarrow 0$; that is, $D_\ff(x^t,x^\star)\rightarrow 0$ in quadratic mean. 

The Bregman divergence is convex in its first argument, so that for every $T\geq 0$, 
\begin{equation*}
D_\ff(\bar{x}^T,x^\star) \leq \frac{1}{T+1} \sum_{t=0}^T D_\ff(x^t,x^\star).
\end{equation*}
Combining this last inequality with \eqref{eqjghhrg} yields 
\begin{equation*}
(T+1)(2-\gamma L_\ff)\Exp{D_\ff(\bar{x}^T,x^\star)} \leq \Psi^0.
\end{equation*}

Now, if $\mu_{\hc}>0$, then $c<1$ in \eqref{eqsm2}, and since $\Psi^t$ converges almost surely to $\Psi^\infty$, 
it must be that $\Exp{\sqnorm{u^{t}-u^\star}}\rightarrow 0$.\hfill$\square$\bigskip

The counterpart of Theorem 2 in the convex case is:\medskip

\noindent\textbf{Theorem 12.}\ \  {\itshape Suppose that $\g=0$, and that $\lambda_{\min}(KK^*)>0$ or $\mu_{\hc}>0$. In \algn{RandProx}, suppose that $0<\gamma< \frac{2}{L_{\ff}}$, $\tau>0$, and $\gamma \tau \big((1-\zeta)\|K\|^2+\oma\big)\leq 1$. 
Then there is a unique dual solution $u^\star$ to \eqref{eqpbd} and $(u^t)_{t\in\mathbb{N}}$ converges to $u^\star$, in quadratic mean.}\medskip

\noindent\textit{Proof of Theorem 12}\ \  Considering the proof of Theorem 2, the same arguments as in the proof of Theorem 11 apply, with $c$ in \eqref{eqsm2} now equal to 
\begin{equation*}
c=1-\frac{2\tau \mu_{\hc} + \gamma\tau\lambda_{\min}(KK^*)}{(1+\omega)(1+2\tau \mu_{\hc})}<1.
\end{equation*}
Hence, $\Exp{\sqnorm{u^{t}-u^\star}}\rightarrow 0$.\hfill$\square$

\end{document}